%% file: Author_tex.tex
\documentclass[12pt]{amsart}

\usepackage{xypic}
\usepackage{amsmath}
\usepackage{amssymb}
\usepackage{enumerate}
\usepackage{color}
\usepackage{pgf}
\usepackage{tikz}
\usetikzlibrary{arrows,automata}
\usepackage{relsize}
\usepackage{latexsym}
\usepackage{pgfplots}
\usetikzlibrary{positioning}
\usepackage{soul} 
\usepackage{amsfonts}
\usepackage[T1]{fontenc}
\usepackage[latin9]{inputenc}
\usepackage{geometry}
\usepackage{amsmath}
\usepackage{amsthm}
\usepackage{stmaryrd}
\usepackage{stackrel}
\usepackage[all]{xy}
\usepackage{xcolor}
\usepackage{graphicx}
\theoremstyle{remark}

\makeatletter
\theoremstyle{plain}
\newtheorem{thm}{\protect\theoremname}[section]
\theoremstyle{definition}
\newtheorem{example}[thm]{\protect\examplename}
\theoremstyle{definition}
\newtheorem{defn}[thm]{\protect\definitionname}
\theoremstyle{plain}
\newtheorem{lem}[thm]{\protect\lemmaname}
\ifx\proof\undefined
\newenvironment{proof}[1][\protect\proofname]{\par
	\normalfont\topsep6\p@\@plus6\p@\relax
	\trivlist
	\itemindent\parindent
	\item[\hskip\labelsep\scshape #1]\ignorespaces
}{%
	\endtrivlist\@endpefalse
}
\providecommand{\proofname}{Proof}
\fi
\theoremstyle{plain}
\newtheorem{prop}[thm]{\protect\propositionname}
\theoremstyle{remark}

\date{}

\makeatother

\usepackage{babel}
\providecommand{\definitionname}{Definition}
\providecommand{\examplename}{Example}
\providecommand{\lemmaname}{Lemma}
\providecommand{\propositionname}{Proposition}
\providecommand{\remarkname}{Remark}
\providecommand{\theoremname}{Theorem}

\begin{document}
\title{Persistent Homology on a lattice of multigraphs}

\author[Joaqu\'in D\'iaz Boils]{Joaqu\'in D\'iaz Boils}
\address{Departament d'Economia Aplicada\\
Facultat d'Economia\\
Avinguda dels Tarongers\\
Universitat de Val\`encia\\
46022-Val\`encia. Spain.}
\email{joaquin.diaz@uv.es}



\maketitle

\begin{abstract}
{A multicomplex structure is defined from an ordered lattice of multigraphs. This structure will help us to observe the features of Persistent Homology in this context, its interaction with the ordering and the repercussions of the process of merging multigraphs in the calculation of the Betti numbers. For the latter, an extended version of the \emph{incremental algorithm} is provided.

The ideas here developed are mainly oriented to the original example described in \cite{Sign} and used more extensively in \cite{Sign2} in the context of the formalization of the notion of \emph{embodiment} in Neuroscience.}

\\
2020 Math Subject Classification: Primary 62R40, Secondary 55U10.
\end{abstract}

\section{Introduction}\label{intro}

This paper studies the mathematical structure introduced in \cite{Sign} and used in a deep way in \cite{Sign2} to provide an algebraic framework for an abstract notion of \emph{embodiment} in Neuroscience by means of multigraphs. Its aim is to examine the presence and extent of Homology in a certain multilayer structure inherited from multigraphs, that is, graphs for which we allow to have parallel edges between the same vertices. For that purpose, we need to translate many of the concepts defined in the references above to the language of \emph{simplicial complexes}. In this line, we hope to contribute showing that the concepts of \emph{Topological Data Analysis} help to give clarity to the formal models of Neuroscience.

The structure considered entails dynamical behavior which needs a suitable setting to be studied. In the present paper, we try to understand which geometric features are preserved along that dynamism. For that purpose, a way to combine multigraphs is defined by obtaining larger and larger networks, that is, the role of the operator $\odot$. Our idea here is to investigate how much the complex structure changes due to the process of merging graphs, in particular to try to know to what extent some features are preserved as long as the dimension of the complex grows through several applications of $\odot$. 

As a general idea, we obviate the classical interaction of an static multilayer network containing graphs (as in \cite{Kivela}). Instead, we focus on its underlying algebraic structure: how it changes the whole network and what remains the same in it. It is precisely through the interaction among graphs using $\odot$ that we obtain a sequence of configurations that allows us to observe these remaining features. 

These sequences are called \emph{filtrations} in the context of a simplicial complex (see, for example, \cite{Akt}). One of the main contributions of this work is to make interact simplices rather than graphs. That is, we construct a simplicial complexes out of graphs and make $\odot$ to act on the simplices obtained by the graphs in the multilayer network. However, rather than considering the usual concept of filtration, we introduce the concept which we call \emph{interaction filtration}, which is nothing more than the sequence of configurations obtained by successive use of our operator $\odot$. 

It is precisely to the aim of obtaining new applications in Neuroscience, that new mathematical ideas are developed which fit with our purpose. In this paper we introduce the construction of multicomplexes out of multigraphs. Multicomplexes were nicely defined in \cite{Lubot}, where many algebraic considerations were made starting from simple graphs rather than multigraphs. Making use of multigraphs, however, we achieve to endow the multicomplexes with much interesting structure and extend their potential applications. 

Moreover, the graphs being here used have the particularity that they are coloured on edges. This is a particularity just mentioned in \cite{Lubot} but necessary for the purpose mentioned of formalizing the structure of networks for which the colours relate to different layers. All this features (multigraphs plus colours) make the multicomplexes richer in their potentiality to generate useful formalizations for concrete examples.

Once the structure of simplicial complex is presented, we investigate usual concepts such as Persistent Homology and the \emph{Betti numbers} to figure out the shape of the network through time. These measures, belonging to the field of Topological Data Analysis, have not been investigated in the context of multicomplexes, to the knowledge of the author, and are addressed in the last pages of the paper. Our contribution here consists of an extension of the \emph{incremental algorithm} given in \cite{Incremental}. The formulae found in that reference, allowing the calculation of the Betti numbers of a series of simplicial configurations in a filtration, are extended in Proposition \ref{formulae} to cover the case of multicomplexes.

As a general purpose, we introduce along the paper as many pictures as possible to facilitate the understanding of the content. On the other hand, the concepts used in the second section, mainly belonging to the Lattice Theory, are to be seen as a mathematical foundation of the structure developed further in the following sections. 

We outline the paper as follows. In Section 2 we introduce the main concepts and present the inspiring example. The crucial idea of ordering is defined and some results emerging from the lattice structure it gives rise to are proved. Section 3 is devoted to deeply explain how multicomplexes can be defined out of multigraphs and how are they merged through the operator $\odot$ working now on simplices. Some categorical considerations are made for the structure (\cite{Baezetal}). A new idea of filtration, the \emph{interaction filtration}, is defined in Section 4 closely related to the graphs obtained by the application of $\odot$. This suggests the use of Homology in the following section by introducing the idea of many-dimensional hole. Section 5 contains a way to calculate Betti numbers for the structure defined as an extension of the known as \emph{Incremental method}, an algorithm defined in \cite{Akt}.  

\section{An ordered set of multigraphs}
    
Let $X$ be a set, a multisubset is a pair $(Y,m)$ where $Y$ is the underlying 
subset of $X$ and $m\colon Y\to\mathbb{Z}^+$ is the \emph{multiplicity function} that 
assigns to each element in $Y$ the number of occurrences (see \cite{Bliz}). 
Next definitions can be found in \cite{Sign} which are based in \cite{Baezetal} 
and \cite{McLane}, where multilayer networks are included into an abstraction 
called \emph{network model}.

\begin{defn}A \emph{multigraph} $G$ on a set of \emph{nodes} $V(G)$ is a multisubset 
of \emph{edges} $E(G)$ that corresponds to pairs of elements of $V(G)$, 
together with the multiplicity function $m_G\colon E(G)\rightarrow\mathbb{Z}^+$. 
\end{defn}
Similarly, the edges could have different colors. Let $C$ be a finite set of colors, 
and $col_G\colon E(G)\to \mathcal{P}(C)$ a mapping that assigns to each edge a 
subset of colors. We will consider a coloured multigraph as the pair $(G,col_G)$ and, for $k\in\mathbb{N}$, 
we say that a multigraph is \emph{$k$-colored} if $col_G$ is onto and $k=|C|$, i.e.
$k$ denotes the number of colors included into the multigraph. We define the set of nodes indexed by the set $\{1,\dots,n\}$ and denote by $\mathcal{G}(n)$ the set of coloured multigraphs with such $n$ nodes. Let $c$ be a single colour, then we denote by $\mathcal{G}^{c}(n)$ the set of $1$-colored layers. Every multigraph is from now on a coloured multigraph.
\subsection{Chains of graphs out of $\otimes$ and $\odot$}
Let $C=\{c_1,\dots,c_k\}$ be a set of colours, then we define the set of multilayer networks as the product

$$
\mathcal{G}^{\otimes C}(n):=(\mathcal{G}^{c_1}\otimes\cdots\otimes \mathcal{G}^{c_k})(n)
= \mathcal{G}^{c_1}(n)\times\cdots\times \mathcal{G}^{c_k}(n)
$$
Then every multigraph in $\mathcal{G}^{\otimes C}(n)$ is called a \emph{$|C|$-coloured multilayer}. We can also have sets of multilayers in the form $\mathcal{G}^{\otimes C}(n) \times \mathcal{G}^{\otimes C'}(m)$ as shown in the following example.

\begin{example}\label{example2}
    A multilayer $G \otimes H$ built out of layers $G \in (\mathcal{G}^{blue} \otimes \mathcal{G}^{red})(3)$ and $H \in (\mathcal{G}^{green} \otimes \mathcal{G}^{yellow})(4)$ can be represented as

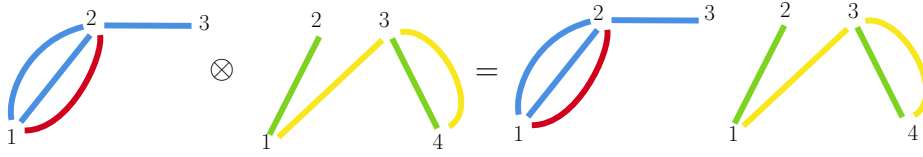
\tiny{
\input{exampletensor.tikz}
}

\end{example}

The tensor product $\otimes$ places layers in parallel establishing a determined order (i.e. $G\otimes H$ is not the same as $H\otimes G$). This operation is non-commutative and imposes a first multilayer order configuration. To simplify notation, we fix nodes into the finite set $V$ and denote such a set by $\mathcal{G}^{\otimes C}$ (without reference to $n$).

Now we are in position to define a commutative binary operation in $\mathcal{G}^{\otimes C}$. Let us denote by $\sqcup$ the disjoint union of sets.

\begin{defn}\label{odot}
Let $C$ and $V$ be fixed sets of colors and nodes respectively.
Let a $k$-colored multigraph $G\in \mathcal{G}^{\otimes C}$ and a $j$-colored multigraph 
$H\in \mathcal{G}^{\otimes C}$, and assume that $C_1:=col(E(G))\subseteq C$, 
$C_2:=col(E(H))\subseteq C$ and that $V(G),V(H)\subseteq V$. Then the operation
\[
\odot\colon \mathcal{G}^{\otimes C}(n)\times \mathcal{G}^{\otimes C}(m)\longrightarrow \mathcal{G}^{\otimes C}
\] 
produces a new ($k+j-s$)-colored multigraph $G \odot H$, where $s=|C_1\cap C_2|$ with $n+m-p$ vertices 
where $p=\left|V(G)\cap V(H)\right|$ defined as 
$V(G\odot H):=V(G)\sqcup V(H)$, $E(G \odot H):=E(G) \cup E(H)$,
$m_{G \odot H}:=m_G + m_H$ and $col_{G \odot H}:=col_G \cup col_H$, where
the mappings are defined in a natural way.
\end{defn}

\begin{example}
For $n=3,m=4,k=j=2,s=0$ and $p=3$:
\vspace{1em}
\tiny{
\input{odot.tikz}
}
\end{example}

It is seen that $\odot$ is a commutative operation while $\otimes$ is not. We establish $\odot$ having priority over $\otimes$,
that is: 
\[
G\otimes H\odot K=G\otimes(H\odot K)
\]

Notice we have defined two different ways of composing multigraphs: $\otimes$ and $\odot$. Also notice that, with this notation, the tensor means no interaction between multigraphs, they are just interpreted as \emph{put together}. That is, we have sets of concatenations in the form 
$$G_1\oslash^{1}\cdots\oslash^{k-1}G_k$$ with $\oslash^{i}\in\{\otimes,\odot\}$ 
for $i=1,\dots,k-1$. 

By abusing the notation we denote by $\mathcal{G}^{\otimes C}(k)$ the set of all chains in the above form constructed out of a fixed set of monochrome multigraphs $G_1,...,G_k$. 

\begin{example}\label{ex1}
For $k=3$ we have the concatenations 
\small{
\[
\begin{array}{c}
\mathcal{G}^{\otimes C}(3):=\{G\otimes H\otimes K,G\otimes K\otimes H,H\otimes G\otimes K,H\otimes K\otimes G,K\otimes G\otimes H,K\otimes H\otimes G,\\
G\odot H\otimes K,G\odot K\otimes H,H\odot K\otimes G,G\otimes H\odot K,H\otimes G\odot K,K\otimes G\odot H,G\odot H\odot K\}
\end{array}
\]}
\end{example}

However, and in order to simplify 
the notation, we sometimes will use lowercase letters as multigraphs and fix the number of colours as $k$, that is: $\mathcal{G}^{\otimes C}(k):=\mathcal{G}^{\otimes C}$ when no other parameter is stated.

\subsection{The partial ordered structure}\label{pomonoid}

The operation $\odot$ defined can be seen as 
an \emph{accumulation of vertices and edges} of two given multigraphs. By using $\odot$ we obtain a new multigraph with possibly more colors than the original ones.

Given $G\otimes H\otimes K$ and $G\odot H\otimes K\in \mathcal{G}^{\otimes C}$
we understand that $G\odot H\otimes K$ is more complex or \emph{is over} $G\otimes H\otimes K$. Following this intuition, we say by convention that $G\otimes H\otimes K\le G\odot H\otimes K$ since
we consider that $G\odot H$ is more complex, in some sense, 
than $G\otimes H$. Let us formalize this idea.
\begin{defn} (from \cite{Birk}) A partially ordered set or a \emph{poset} is a set with a binary operation
$\le$ which is reflexive, antisymmetric and transitive. Let $P$ be a poset. We say that $b \in P$ is a \emph{bottom} element if 
$b\le x$ for every $x \in P$ and $a\in P$ is a \emph{top} element
if $a\ge x$ for every $x\in P$.
\end{defn}
Now we define the relation $\le$ in $\mathcal{G}^{\otimes C}$ by ordering the concatenations of 
multigraphs as given in the following.

\begin{defn}\label{deforder}
Let $\pi$ be a permutation of the set $\{1,\dots,k\}$ with $k=|C|$. For every $\oslash^{l} , \ominus^{l} \in \{\otimes,\odot\}$ with $1\leq l \leq k-1$ we write
\[
G_{\pi(1)}\oslash^{1}\cdots\oslash^{k-1} G_{\pi(k)}\le G_{\pi(1)}\ominus^{1}\cdots\ominus^{k-1} G_{\pi(k)}
\]
if and only if there is no $l$ such that $\oslash^{l}=\odot$
and $\ominus^{l}=\otimes$.

\end{defn}

\begin{example}\label{example_chains}
The partial character of this ordering is shown by observing that:
$$G\otimes H \otimes K \leq G\otimes H \otimes K \otimes L \leq G\odot H \otimes K \otimes L$$
and
$$G\otimes H \otimes K \leq G\odot H \otimes K \leq G\odot H \otimes K \otimes L$$
while
$$G\odot H \otimes K \not\leq G\otimes H \otimes K \otimes L$$ 
and
$$G\otimes H \otimes K \otimes L \not\leq G\odot H \otimes K$$

\vspace{1em}
\end{example}

This partial order allows us to define the following mappings. 

\begin{defn}
Let the mapping $f_{j}\colon \mathcal{G}^{\otimes C}\to \mathcal{G}^{\otimes C}$:
\[
f_{j}(x)=\begin{cases}
G_{1}\oslash^{1}\cdots G_{j}\odot G_{j+1}\cdots\oslash^{k-1}G_{k} & 
\textrm{if }x=G_{1}\oslash^{1}\cdots G_{j}\otimes G_{j+1}\cdots\oslash^{k-1}G_{k}\\
x & \textrm{otherwise}
\end{cases}
\]
for $j=1,\dots,k-1$. We say that 
$x,y\in \mathcal{G}^{\otimes C}$ are \emph{comparable through $f_j$} if $f_j(x)=y$.
\end{defn}

By adding $f_0$ as the identity, it is easy to see that $f_{j}$ are order-preserving.
For the sake of clarity we use the notation $f_{j}$ for any mapping defined above,
avoiding the list of indices. These mappings will be useful in the sequel, the next example illustrates how these functions 
work and describe, in some sense, a flow on $\mathcal{G}^{\otimes C}$ as a poset.

\begin{example}\label{ex2}
For the elements in Example \ref{ex1} we have: 
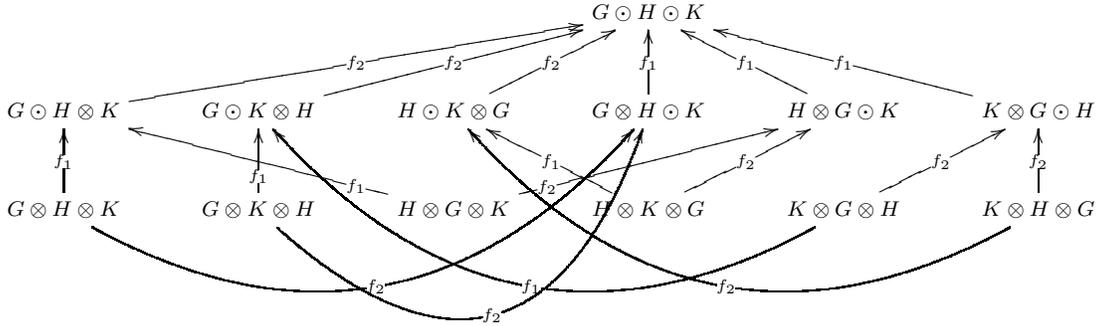
\tiny{
{
\[
\xymatrix{ &  & & G\odot H\odot K\\
G\odot H\otimes K\ar[urrr]|{f_{2}} & G\odot K\otimes H\ar[urr]|{f_{2}} & 
H\odot K\otimes G\ar[ur]|{f_{2}} & G\otimes H\odot K\ar[u]|{f_{1}} & 
H\otimes G\odot K\ar[ul]|{f_{1}} & K\otimes G\odot H\ar[ull]|{f_{1}}\\
G\otimes H\otimes K\ar[u]|{f_{1}}\ar@/_{4pc}/[urrr]|{f_{2}} & 
G\otimes K\otimes H\ar[u]|(.35){f_{1}}\ar@/_{5pc}/[urr]|{f_{2}} & 
H\otimes G\otimes K\ar[ull]|(.25){f_{1}}\ar[urr]|(.24){f_{2}} & 
H\otimes K\otimes G\ar[ul]|(.50){f_{1}}\ar[ur]|{f_{2}} & 
K\otimes G\otimes H\ar@/^{4pc}/[ulll]|{f_{1}}\ar[ur]|{f_{2}} & 
K\otimes H\otimes G\ar@/^{4pc}/[ulll]|{f_{2}}\ar[u]|{f_{2}}
}
\]
}
}

\end{example}

From the example above we extract two immediate results.
The first one establishes that one can obtain the top element after
an action of every $f_{j}$ over a given concatenation whatever ordering could
be and the second that $f_{j}$ are increasing.

\begin{lem}
$f_{i_{1}}\cdots f_{i_{k}}(G_{1}\oslash^{1}\cdots\oslash^{k-1}G_{k})=G_{1}\odot\cdots\odot G_{k}$
for $i_{1}<\cdots<i_{k}$ a permutation of $1,\dots,k$.
\end{lem}
\begin{lem}\label{monoto}
$f_{j}(G_{1}\oslash^{1}\cdots\oslash^{k-1}G_{k})\geq G_{1}\oslash^{1}\cdots\oslash^{k-1}G_{k}$.
\end{lem}


 
\begin{prop}\label{prop1}
$\mathcal{G}^{\otimes C}$ is a partial ordered set with top $G_{1}\odot\cdots\odot G_{k}$.
\end{prop}
\begin{proof}
Observe the order of $\mathcal{G}^{\otimes C}$ is described by the mappings $f_j$ 
(see Example \ref{ex2}). Reflexivity is given by $f_0$ while transitivity is 
immediate by definition of the mappings $f_j$. For antisymmetry 
we recall the form of the ordering given in the previous definition, 
now a concatenation can only be compared both ways with another concatenation if they 
are both the same. In that case they are compared by means of the 
same $f_j$ whenever a $\odot$ appears in the $j$-position of the concatenation. 
\end{proof}

Notice that we cannot dualize the above since inverse mappings in  such as $g_{1}$ for which 
$$g_{1}(G\odot H\otimes K)=G\otimes H\otimes K$$
lose the well-definedness condition for the non commutativity of $\otimes$.\\


\subsection{A lattice-ordered partial monoid}\label{lattice}

We saw above that we can define a partial order in a set of multigraphs
and showed that this order yields several properties in such a framework.
In this section we go a little further and provide a lattice structure for
$\mathcal{G}^{\otimes C}$. 

A meet (resp. join) semilattice is a poset $(L,\le)$ such
that any two elements $x$ and $y$ have a greatest lower bound (called
meet or infimum) (resp. a smallest upper bound (called join or supremum)),
denoted by $x\wedge y$ (resp. $x\vee y$).
A poset $(L,\le)$ is called a \emph{lattice} and denoted by 
$(L,\le,\wedge,\vee)$ if for every pair of elements we can construct 
into the lattice their meet and their join. These definitions can be 
found for instance in \cite{Birk}.
Let us define a meet and a join operators for the poset $(\mathcal{G}^{\otimes C},\leq)$:
\begin{defn}
Given $G_{1}\oslash^{1}\cdots\oslash^{k-1}G_{k}$ and  
$G_{1}\ominus^{1}\cdots\ominus^{k-1}G_{k}\in \mathcal{G}^{\otimes C}$
(in short $\oslash G$ and $\ominus G$) we write $\oslash G\wedge\ominus G=\olessthan G$
for $G_{1}\olessthan^{1}\cdots\olessthan^{k-1}G_{k}$ such that 
\[
\olessthan^{j}=\begin{cases}
\otimes & \textrm{if }\oslash^{j}=\otimes\textrm{ or }\ominus^{j}=\otimes\\
\odot & \textrm{otherwise}
\end{cases}
\]
and we write $\oslash G\vee\ominus G=\ogreaterthan G$ for $G_{1}\ogreaterthan^{1}\cdots\ogreaterthan^{k-1}G_{k}$
such that 
\[
\ogreaterthan^{j}=\begin{cases}
\odot & \textrm{if }\oslash^{j}=\odot\textrm{ or }\ominus^{j}=\odot\\
\otimes & \textrm{otherwise}
\end{cases}
\]
\end{defn}

One can easily check the usual properties of both operations. That is, for every $x,y,z\in \mathcal{G}^{\otimes C}$: $x\wedge y\leq x,y$ and for every $z\leq x,y$
one has $z\leq x\wedge y$. And dually: $x\vee y\geq x,y$ and for
every $z\geq x,y$ one has $z\geq x\vee y$.

\begin{prop}\label{absorbe}
The \emph{absorption laws} are satisfied for every $x,y\in \mathcal{G}^{\otimes C}$:
\begin{itemize}

\item{$x\vee (x\wedge{y})=x$}
\item{$x\wedge (x\vee{y})=x$}
\end{itemize}

\end{prop}
\begin{proof}
Let us prove the first assertion. For $x=\ominus G$ and $y=\obar G$ we construct $x\wedge{y}=\oslash G$ such that 
\[
\oslash^{j}=\begin{cases}
\otimes & \textrm{if }\obar^{j}=\otimes\textrm{ or }\ominus^{j}=\otimes\\
\odot & \textrm{otherwise}
\end{cases}
\] and $x\vee (x\wedge{y})=\ogreaterthan G$ as \[
\ogreaterthan^{j}=\begin{cases}
\odot & \textrm{if }\ominus^{j}=\odot\textrm{ or (}\obar^{j}=\odot\textrm{ and }\ominus^{j}=\odot)\\
\otimes & \textrm{otherwise}
\end{cases}
\] which can be expressed as 
\[
\begin{cases}
\odot & \textrm{if }\ominus^{j}=\odot\\
\otimes & \textrm{otherwise}
\end{cases}
\]
and becomes the same assignation as that considered for $x=\ominus G$.
\end{proof}

Recalling that a \emph{minimal element} in a 
poset is an element such that it is not greater than any other element in the poset we have:
\begin{prop}
$(\mathcal{G}^{\otimes C},\leq,\wedge,\vee,1_{\odot},m_{\pi})$ is
an upper-bounded lattice where:
\begin{itemize}
\item $1_{\odot}=G_{1}\odot\cdots\odot G_{k}$ is the top and 
\item $m_{\pi}=G_{\pi(1)}\otimes\cdots\otimes G_{\pi(k)}$ are $k!$ 
minimal elements for $\pi$ a permutation of the set $\{1,\dots,k\}$.
\end{itemize}
\end{prop}
\begin{proof}
Check that $x\wedge1_{\odot}=x,x\vee1_{\odot}=1_{\odot}$, $x\wedge m_{\pi}=m_{\pi}$ and $x\vee m_{\pi}=x$.
\end{proof}
\begin{prop}
$(\mathcal{G}^{\otimes C},\leq,\wedge,\vee,1_{\odot},m_{\pi})$ is
distributive.
\end{prop}
\begin{proof}
Let $x_{1}=\oslash G$, $x_{2}=\ominus G$, $x_{3}=\obar G$.
Now $x_{1}\wedge(x_{2}\vee x_{3})=\olessthan G$ where
\[
\olessthan^{j}=\begin{cases}
\otimes & \textrm{if }\ominus^{j}=\otimes\textrm{ or }\obar^{j}=\otimes\textrm{ and }\oslash^{j}=\otimes\\
\odot & \textrm{if }\ominus^{j}=\obar^{j}=\odot\textrm{ or }\oslash^{j}=\odot
\end{cases}
\]
which is exactly the same operator as
\[
\begin{cases}
\otimes & \textrm{if}\textrm{ no }(\oslash^{j}=\odot\textrm{ or }\ominus^{j}=\odot)\textrm{ or }\textrm{no }(\oslash^{j}=\odot\textrm{ or }\obar^{j}=\odot)\\
\odot & \textrm{otherwise}
\end{cases}
\]
for $(x_{1}\wedge x_{2})\vee(x_{1}\wedge x_{3})$.
\end{proof}
\begin{prop}\label{prop2}
Mappings $f_{j}$ preserve meets and joins.
\end{prop}
\begin{proof}
Straightforward calculations.
\end{proof}
\begin{subsection}{Complements} 
In \cite{Birk} a \emph{complemented lattice} is defined as a 
bounded lattice (with least element 0 and greatest element 1)
in which every element $a$ has a \emph{complement}, i.e. an element 
$b$ such that $a \vee {b} = 1$ and $a \wedge {b} = 0$.
Also, given a lattice $L$ and $x\in{L}$ we say that $\hat{x}$ is 
an \emph{orthocomplement of $x$} if the following conditions are satisfied:
\begin{itemize}
\item{$\hat{x}$ is a complement of $x$}
\item{$\hat{\hat{x}}=x$}
\item{if $x\leq{y}$ then $\hat{y}\leq{\hat{x}}$}.
\end{itemize}
A lattice is \emph{orthocomplemented} if every element has an orthocomplement. We give a slightly different approach:

\begin{defn}
We say that an upper bounded lattice $(L,\leq,\wedge,\vee,1)$ with a set of minimal elements $\{m_1,...,m_k\}$ is \emph{semi-orthocomplemented} if every element $a\in {L}$ has a complement, i.e. an element 
$b$ such that $a \vee {b} = 1$ and $a \wedge {b} = m_i$ for a certain $i\in \{1,...,k\}$.
\end{defn}

\begin{prop}
$(\mathcal{G}^{\otimes C},\leq,\wedge,\vee,1_{\odot},m_{\pi})$ is a semi-orthocomplemented lattice.
\end{prop}
\begin{proof}
For $x=\oslash G$ consider $\hat{x}=\ominus G$ where
\[
\ominus^{j}=\begin{cases}
\otimes & \textrm{if }\oslash^{j}=\odot\\
\odot & \textrm{if }\oslash^{j}=\otimes\,,
\end{cases}
\]
\end{proof}
\end{subsection}
\subsection{Partiality}
  Now some concepts from \cite{Jasem} are taken and adapted for the case of a partial operation such as ours (since not every pair of elements are comparable through $\leq$). 

\begin{defn}
    A system $(A,+,\leq,\wedge,\vee)$ is called a \emph{lattice-ordered partial monoid} if 
    \begin{itemize}
        \item{$(A,+)$ is a partial monoid}
        \item{$(A,\leq)$ is a lattice with $\wedge$ and $\vee$}
        \item{$a\leq{b}$ implies $a+x\leq{b+x}$ and $x+a\leq{x+b}$}
        \item{$a+(b\vee{c})=(a+b)\vee{(a+c)}, (b\vee{c})+a=(b+a)\vee{(c+a)}$} 
        \item{$a+(b\wedge{c})=(a+b)\wedge{(a+c)}, (b\wedge{c})+a=(b+a)\wedge{(c+a)}$}
    \end{itemize}
for every $a,b,c,x\in{A}$. 
\end{defn}

We are introducing a different operation from that considered in \cite{Jasem}: $+$ will be partial for $\mathcal{G}^{\otimes C}$. This is the reason to study 
it as a lattice-ordered partial monoid. Firstly, a partial semigroup structure will be defined for $\mathcal{G}^{\otimes C}$ by defining $+$ as:

\[
x+y=\begin{cases}
y & \textrm{if } x\geq{y}\\
x & \textrm{if }y\geq{x}
\end{cases}
\]
for $x, y \in \mathcal{G}^{\otimes C}$.

Now we see that $+$ is an associative, commutative and partial operation. It is actually a \emph{partial minimum}.

\begin{prop}
$(\mathcal{G}^{\otimes C},+)$ is a partial commutative monoid.
\end{prop}

\begin{proof}
The operation $+$ satisfies associativity: suppose that 
$x,y,z\in \mathcal{G}^{\otimes C}$ are comparable to each other. Now: $$x+y+z=min(x,y,z)=min(min(x,y),z)=min(x,min(y,z))\,.$$
As $\mathcal{G}^{\otimes C}$ is finite, the unique top element (see Proposition \ref{prop1})
is the identity element.
\end{proof}

Let $(M,\cdot)$ be a partial monoid and let $f\colon M\to M$ 
be a mapping. Recall that a \emph{partial homomorphism} between partial 
monoids is mapping that preserves the binary operation, namely 
$f(x+y)=f(x)+f(y)$, $f(1)=1$ and $x+y\in M$ implies that $f(x)+f(y)\in M$.
A mapping between lattice-ordered partial monoids is a \emph{lattice partial 
homomorphism} if it is a partial homomorphism of partial monoids that
preserves meets and joins.
\begin{prop}
Mappings $f_j$ are partial homomorphisms.
\end{prop}
\begin{proof}
By virtue of Proposition \ref{prop2}, mappings $f_j$ preserve meets and joins. From the definition of the partial operation $+$, we know that $x$ and $y$ are comparable if and only if there exists $x+y$. As $f_j$ is monotone, if $x\le y$, then $f_j(x)\le f_j(y)$
and $f_j(x)$ and $f_j(y)$ are comparable. So 
$f_j(x+y)=f_j(\min(x,y))=f_j(x) = \min(f_j(x),f_j(y))= f_j(x)+f_j(y)$.
Finally, as $1$ is the top element $x\le 1$ for every $x$, $f_j(1)\le 1$. But for Lemma \ref{monoto} we have $1\le f_j(1)$. Therefore $f_j(1)=1$.
\end{proof}

\begin{prop}
$(\mathcal{G}^{\otimes C},+,\leq,\wedge,\vee)$ is a lattice-ordered partial monoid.
\end{prop}

\begin{proof}
    Suppose that $x,y,z\in \mathcal{G}^{\otimes C}$ are comparable. Observe that $$x+(y\vee{z})=(x+y)\vee{(x+z)}, (y\vee{z})+x=(y+x)\vee{(z+x)}$$ and $$x+(y\wedge{z})=(x+y)\wedge{(x+z)}, (y\wedge{z})+x=(y+x)\wedge{(z+x)}$$
together with the fact that for $x\leq{y}$: 
$$
x+z\leq{y+z}, z+x\leq{z+y}\,.
$$
Notice in particular that
$$
x+(y\vee{z})=min(x,y\vee{z})=\begin{cases}
min(x,\odot) & \textrm{if } y=\odot \text{ or } z=\odot\\
min(x,\otimes) & \textrm{else }
\end{cases}
=\begin{cases}
x & \textrm{if } y=\odot \text{ or } z=\odot\\
\otimes & \textrm{else }
\end{cases}
$$
equals to
$$
(x+y)\vee{(x+z)}=min(x,y)\vee{min(x,z)}=\begin{cases}
\odot & \textrm{if } min(x,y)=\odot \text{ or } min(x,z)=\odot\\
\otimes & \textrm{else }
\end{cases}
$$
$$
=\begin{cases}
\odot & \textrm{if } x=y=\odot \text{ or } x=z=\odot\\
\otimes & \textrm{else}
\end{cases}
$$
\end{proof}

\section{Multicomplexes out of multigraphs}

In this section we show how to define multicomplexes out of multigraphs. Our main references for this will be \cite{Akt,Lubot}. The concept of multicomplex is developed in \cite{Lubot} for simple graphs, where the complexes are in addition ordered and coloured in their vertices. We will follow the notation and formulation of that paper but will make use of the concepts introduced so far, the operation $\odot$ in particular, and will colour the edges rather than the vertices.  

\subsection{The clique complex of a graph}
Let us begin this section by sketching some concepts. A \emph{simplicial complex} $\mathcal{K}$ with vertex set $V$ is a non-empty collection of finite subsets of $V$, called \emph{cells} such that it is closed under inclusion. That is: for $A\in \mathcal{K}$ and $B\subseteq A$ then $B\in \mathcal{K}$. The dimension of
a cell $A$ is $|A|-1$ and $\mathcal{K}_j$ denotes the set of \emph{j-cells} (cells of dimension $j$) for $j\geq-1$. We write them as $\{v_0,v_1,...,v_j\}$. The $0-cells$ of a simplicial complex are the vertices themselves.
The dimension of a simplicial complex is the maximal
dimension of a cell in it. For a (j + 1)-cell $\tau = \{\tau_0, . . . , \tau_{j+1}\}$, its boundary $\partial\tau$ is defined to be the set of j-cells $\{\tau - \{\tau_i\}\}^{j+1}_{i=0}$. 

With all these concepts, taken from \cite{Lubot}, we construct a simplicial complex in a different way. Recall that a \emph{k-clique} into a graph $G$ is a complete subgraph of dimension $\mathcal{K}$.

\begin{defn}\label{cliquecomplex} 
Given a graph $G$ its \emph{clique complex}, denoted by $Cl(G)$, has: 
\begin{itemize}
    \item{the vertices of $G$ as its vertices and  }
    \item{the $k-cliques$ as its $(k-1)-cells$.}
\end{itemize} 
\end{defn}

\begin{example}
The following picture:

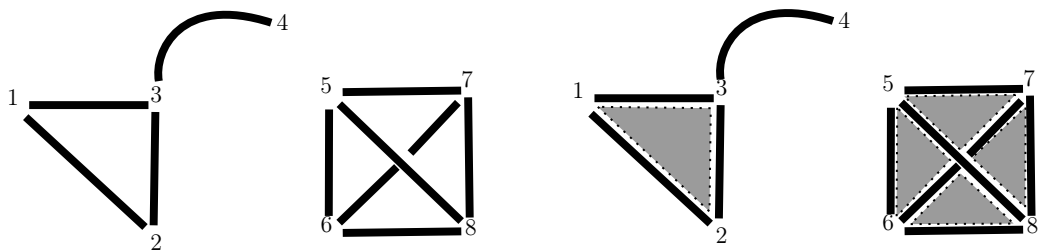
\input{Firstclique.tikz}

\hspace{-1em}shows at the left hand side a 4-clique with set of vertices $\{5,6,7,8\}$, a 3-clique with vertices $\{1,2,3\}$ and a 2-clique with vertices $\{3,4\}$. Then it gives rise (at right) to a complex with a 3-cell, a 2-cell and a 1-cell respectively.

That is, whenever a complete graph is found into one of the connected components of a graph, it converts into a new cell of one less dimension than the clique. In this case we say that a (sub)graph \emph{has been closed} or that a cell \emph{has been created}.

In our notation for the remainder of the paper, we will fill the cells as in the previous example. Observe we can have a picture such as the following, where the square at right has not been closed:
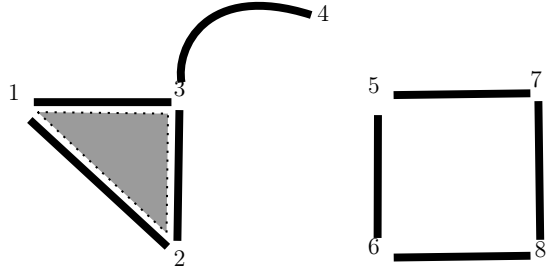
\input{Secondclique.tikz}
\end{example}

\subsection{Multicomplexes out of multigraphs}
Since our goal is to give the structure of a simplicial complex to a set of multigraphs (rather than simple graphs) every definition has to be adapted to that setting. 

Let $\mathcal{K}$ be a simplicial complex with set of vertices $V$ and $m : \mathcal{K} \rightarrow \mathbf{N}$ a \emph{multiplicity function}. Now consider as in \cite{Lubot} the sets $$\mathcal{K}_{m} = \{(\tau, n) / \tau \in \mathcal{K} \text{ and } 1\leq n \leq m(\tau)\}$$ together with the \emph{gluing functions}
$\phi:\{((\tau, n),\sigma) \in \mathcal{K}_m \times \mathcal{K} / \sigma \in \partial \tau\} \rightarrow \mathcal{K}_m$. Observe that $$\phi((\tau, n),\sigma) \in \{(\tau, n) / n \in \{m(\tau)\}\}$$

Here $\phi$ will allow us to identify every copy of a given cell which belongs to a certain boundary of a higher-level cell. Now we are ready to introduce the concept of multicomplex.

\begin{defn}
Considering the definitions above, a 3-tuple $(\mathcal{K},m,\phi)$ for which every $(\tau, n) \in S_m$ and $(\sigma, p), (\rho, q) \in (\tau, n)$ with the same dimension such that $\delta=\sigma \cap \rho$ satisfies $$\phi((\sigma, p), \delta)=\phi((\rho, q),\delta)$$ is a \emph{multicomplex}.
The elements of $\mathcal{K}_m$ are the \emph{multicells} of $(\mathcal{K},m,\phi)$.
  
\end{defn}

By an abuse of notation, we denote simply by $\mathcal{K}$ a multicomplex.

\begin{defn}\label{cliquesimplex}
Cells with the same vertices into a multicomplex are called \emph{duplications}.
\end{defn}
The following definition is an adaptation for multigraphs of an analogous definition in \cite{Akt} for simple graphs. We are from this point deliberately confusing \emph{graph} and \emph{multigraph}.
\begin{defn}\label{cliquecomplex} 
Given a multigraph $(G,col_G)$ its \emph{clique multicomplex}, denoted by $\mathcal{K}_G$, has: 
\begin{itemize}
    \item{the vertices of $G$ as its vertices and  }
    \item{the $k-cliques$ as its $(k-1)-multicells$.}
\end{itemize} 
\end{defn}

From this point we will be working with the clique multicomplex that can be defined in $\mathcal{G}^{\otimes C}$, that is, every clique appearing in a  multigraph into the chains of $\mathcal{G}^{\otimes C}$ is seen now as a multicell. For a certain $x \in \mathcal{G}^{\otimes C}$ we will denote the multicomplex defined this way as $\mathcal{K}_x^{C}$ but we omit the $x$ as often as it does not add any useful information. 

\subsection{Merging complexes}

Once we have the structure of multicomplexes from multi- graphs, we want to define the operation $\odot$ for merging simplices rather than graphs. By abusing the notation we define a binary operation $\odot$ for multicomplexes treated as graphs where we omit the symbols $\otimes$ in the pictures shown. The reason to suppress the tensors in the examples is that, once the simplicial complex is constructed out of an element of $\mathcal{G}^{\otimes C}$, the simplices appear just one next to the other. In other words, when we generate simplicial complexes the initial ordering of the multilayer network is forgotten. 

From the content above, we notice that closing a multigraph has now a different behaviour: it can happen by means of the operation $\odot$ as seen in the following.

\begin{example} \label{basicexample}
In the picture, two 2-simplices have been created (appearing in it with a different tone of grey):
    \input{Thirdclique.tikz}

\hspace{-1em}whose set of vertices is both the same: $\{1,2,3\}$.

Of course we now need to specify each of the 2-cells created in this way since the boundaries for both are the same. That is, the 1-cells $\{1,2\},\{2,3\}$ and $\{3,1\}$. It is precisely the function $\phi$ that allows to identify them: we have now $(\{1,2,3\},1)$ and $(\{1,2,3\},2)$ with boundaries $$(\{1,2\},1),(\{1,3\},1),(\{2,3\},1)$$ and $$(\{1,2\},1),(\{1,3\},2),(\{2,3\},1)$$ respectively.

\end{example}

\subsection{Multiboundaries of a multicomplex}
At this point we note that not every interaction through $\odot$ gives rise to a new cell since it does not always close a graph. On the other hand, a single action of $\odot$ could close more than one cell at a time as happens in the following example.

\begin{example} \label{3cell}
In the following picture, four 2-cells and one 3-cell have been created:

    \input{Forthclique.tikz}
\end{example}

Now for the previous example the following 2-cells have been created, which we detail thanks to $\phi$:

$$\tau_1=(\{1,2,3\},1) \text{ with }\left \lbrace\begin{array}{c}
\phi(\tau_1,\{2,3\})=(\{2,3\},1)\\\phi(\tau_1,\{1,3\})=(\{1,3\},1)\\\phi(\tau_1,\{1,3\})=(\{1,2\},1)
\end{array}\right.$$

$$\tau_2=(\{1,2,3\},2) \text{ with }\left \lbrace\begin{array}{c}
\phi(\tau_2,\{2,3\})=(\{2,3\},2)\\\phi(\tau_2,\{1,3\})=(\{1,3\},1)\\\phi(\tau_2,\{1,3\})=(\{1,2\},1)
\end{array}\right.$$

$$\sigma=(\{1,3,4\},1) \text{ with }\left \lbrace\begin{array}{c}
\phi(\sigma,\{1,3\})=(\{1,3\},1)\\\phi(\sigma,\{1,4\})=(\{1,4\},1)\\\phi(\sigma,\{3,4\})=(\{3,4\},1)
\end{array}\right.$$

$$\rho=(\{2,3,4\},1) \text{ with }\left \lbrace\begin{array}{c}
\phi(\rho,\{2,3\})=(\{2,3\},1)\\\phi(\rho,\{2,4\})=(\{2,4\},1)\\\phi(\rho,\{3,4\})=(\{3,4\},1)
\end{array}\right.$$ 
as well as the 3-cell $(\{1,2,3,4\},1)$.

It is clear from the picture that we cannot talk about boundaries but rather about \emph{multiboundaries}. 

\begin{defn}
We define the multiboundary of a multicell $(\tau, n)$ as $$\partial ^ {m} (\tau, n)= \{\phi((\tau, n),\sigma) / \sigma \in \partial (\tau)\}$$ This set is set to be empty in case the dimension of $(\tau, n)$ is -1.    
\end{defn}

\subsection{Colouring multicomplexes}

The important characterization for us is of course that of edge-coloured multigraphs, to connect with the example introduced in section 2. Therefore, we will continue using the notation introduced in \cite{Lubot} but considering two key modifications:

\begin{itemize}
\item{rather than the vertices, it will be the edges that will be coloured and}
\item{colours could be repeated along the edges of the multigraph as suggested in \cite{Lubot} (\emph{Remark 3.8}.})
\end{itemize}

\begin{defn}
Given a finite dimensional complex $\mathcal{K}$ and a set of colours $C$ we define a \emph{colouring} of its 1-cells as a function $col : \mathcal{K}^{(1)} \rightarrow C$ where $\mathcal{K}^{(1)}$ is the set of 1-cells in $\mathcal{K}$. We colour a finite dimensional multicomplex $(\mathcal{K},m,\phi)$ by colouring the underlying complex $\mathcal{K}$.      
\end{defn}  

With the notation of the previous definition, one can extend the colouring to
all the multicells of $\mathcal{K}$ by means of the multisets $$col(\tau):=\{col(v) : v \in \tau\}$$ for $\tau \in \mathcal{K}^{(j)}$ (where $\mathcal{K}^{(j)}$ is the set of j-cells in $\mathcal{K}$) by an abuse of notation. Therefore, a j-cell into $\mathcal{K}$ is coloured using as many (possibly repeated) colours as 1-cells it contains.

\begin{example}
The following graph
    
    \input{color.tikz}
    
\hspace{-1em}for which $$\mathcal{K}^{(1)}=\{(\{1,2\},1),(\{1,2\},2),(\{1,2\},3),(\{1,3\},1),(\{2,3\},1)\}$$
is coloured through the following assignations:
$$col(\{1,2,3\},1)=\{red, black^{(2)}\}$$
$$col(\{1,2,3\},2)=\{black^{(3)}\}$$
$$col(\{1,2,3\},3)=\{black^{(3)}\}$$
where the superindices indicate multiplicity into the multisets built from assignations such as $$col(\{1,2\},3)=red$$
\end{example}

Since the process of merging graphs through $\odot$ has to avoid duplicating repeated vertices, there is an ordering for them that has to be preserved along the merging process of graphs. That is, we will enumerate the vertices in such a way every one appears just once. 

Since vertices are seen as $0-cells$, successive $j-cells$ are then numbered from them. However, the colours of the 1-cells into the multicells could be repeated and then we need not just a numbering but also a colour to specify them. 

\begin{defn} (from \cite{Lubot})
We denote by $Perm_X$ the permutation group of the set $X$ and $\delta (\tau)=\{\tau \in \mathcal{K} / \sigma \in \partial^{m}\tau\}$. Let $(\mathcal{K},m,\phi)$ be a $d$-dimensional coloured multicomplex. We say that $p =(p_{\tau})$ is an \emph{ordering} if for $\tau \in \mathcal{K}^{(d-1)}$ the assignation $$p:\mathbf{Z} \rightarrow Perm_{\delta (\tau)}$$ is a group homomorphism such that for every $\sigma, \sigma ' \in \delta (\tau)$ there exists $n \in \mathbf{Z}$ such that $p(n). \sigma=\sigma'$ where $p(n).\sigma$ is the action of the permutation $p(n)$ on the $d-multicell$ $\sigma$.
\end{defn}

Now we need to avoid ambiguities such as the following:

\begin{example}
In the graph
    
    \input{color2.tikz}

\hspace{-1em}we have to distinguish between the 2-cells $\tau_1=(\{1,2,3\},1)$ such that $col(\phi(\tau_1,\{1,2\}))=red$ and $\tau_2=(\{1,2,3\},\{1,2\})$ such that $col(\phi(\tau_2,\{1,2\}))=black$ and for which $$col(\tau_1)=col(\tau_2)=\{red^{(2)}, black\}$$ 
\end{example}


That is the reason we must consider an orientation of the $j-cells$ in a multicomplex. Therefore, from now on the function $col$ is thought to be defined on an ordered ${\mathcal{K}^{(1)}}$ and therefore the sets of colours generated by it will be also ordered without any further reference. For instance, in the previous example: $col(\tau_1)=\{red,black,black\}$ and $col(\tau_2)=\{black,black,red\}$.

\subsection{The monoidal category of multigraphs}
From an algebraic perspective, consi- dering the operation $\odot$ for the simplices rather than the graphs entails the definition of a tensor product for the category of multicomplexes. We have in fact something stronger: it gives to that category the structure of a \emph{strict symmetric monoidal category}.      

The category of coloured multicomplexes was defined in \cite{Lubot}, we adapt it in a suitable manner for our context and with the notation here introduced. 

\begin{defn}
    Let $\mathbb{K}$ the category whose
    \begin{itemize}
    \item{objects are tuples $(\mathcal{K}_x,col,p)$ for $x \in \mathcal{G}^{\otimes C}$, $col$ a colouring and $p$ a permutation}
    \item{morphisms between $(\mathcal{K}_x,col,p)$ and $(\mathcal{K}_y,col',p')$ are simplicial multimaps $\phi$ preserving the structure. That is: 
    
    -$\phi(x)=y$
    
    -$col'\circ \phi=col$
    
    -$\phi(p(n).\sigma)=p'(\phi(n)).\phi(\sigma)$}
    \end{itemize}
\end{defn}

$\mathbb{K}$ is the category of multicomplexes to which $\odot$ gives a well-known structure.

\begin{defn}
    A \emph{strict symmetric monoidal category} $(\mathbb{C},\boxtimes,1)$ is a category $\mathbb{C}$ equipped with an object $1 \in \mathbb{C}$ and a bifunctor $\boxtimes:\mathbb{C}\times \mathbb{C}\rightarrow\mathbb{C}$ such that for objects $A,B,C$ and morphisms $f,g,h$:

    -$(A\boxtimes B)\boxtimes C=A\boxtimes (B\boxtimes C)$
    
    -$1\boxtimes A=A$

     -$A\boxtimes B=B\boxtimes A$
    
    -$(f\boxtimes g)\boxtimes h=f\boxtimes (g\boxtimes h)$

    -$f\boxtimes g=g\boxtimes f$
    
    -$Id_1\boxtimes f=f$

\end{defn}

\begin{lem}
$(\mathbb{K},\odot,\emptyset)$ is a strict symmetric monoidal category.
\end{lem}

\begin{proof}
By definition \ref{odot} we have associativity and commutativity. The empty complex $\emptyset$ plays the role of a neutral element.
\end{proof}

Observe that in this section we considered multicomplexes consisting on a unique connected component, that is, not coming from chains of graphs but rather from a single multigraph. 

\section{Filtrations}
It is known (see \cite{What}) that a \emph{filtration} of a simplicial complex $\mathcal{K}$ is just a nested (and totally ordered) set of simplicial complexes $\mathcal{K}_i$ contained in $\mathcal{K}$. That is, a sequence
$$\emptyset =\mathcal{K}_0 \subseteq \mathcal{K}_1 \subseteq...\subseteq \mathcal{K}$$
Of course, there is a number of ways of ordering the different $\mathcal{K}_i$ (see \cite{Akt} for a nice listing). Filtrations are considered in the context of simplicial complexes as a way to describe and organize the construction of a simplicial complex through a discrete process. 

We formalize this \emph{complexification} process by defining a different concept of filtration. For a certain $\mathcal{K}^{C}_x$ we consider the complexes generated by the successive interactions of multigraphs through $\odot$ starting from $x=\ominus G$. In the following we will omit the sub- and the superscript of $\mathcal{K}$ when the initial graph $x$ is understood from the context or the set of colours $C$ is not relevant.





\subsection{The interaction filtration}
This paper introduces a sort of filtration not having the form of a totally ordered set, as every other filtration defined so far to the knowledge of the author, but that of a partial ordered set. As explained below, this ordering for $\mathcal{K}^{C}$ is inherited by the posetal behaviour of $\mathcal{G}^{\otimes C}$. 


Let us describe a \emph{filtration} for $\mathcal{K}^{C}$, for that purpose we see the chains as elements of the complex. Given $\oslash{G}=G_{1}\oslash^{1}\cdots\oslash^{k-1}G_{k}\in \mathcal{G}^{\otimes C}$ and $[k]=\{1,...,k-1\}$ let be $$\xi=|\{j\in [k] \text{ such that } \oslash^{j}=\odot\}|$$ the numbers of times $\odot$ appears in a chain for the rest of the section. 

We will consider for simplicity that every graph $G_i$ into $x$ has a single connected component. Once the simplicial complex is constructed out of an element into $\mathcal{G}^{\otimes C}$, there is no reason to consider $\otimes$ as a non-commutative operation any more. 

For $j \in[k]$ the subset
$$\mathcal{K}^{C}_x[j]=\{\oslash{G}\in \mathcal{K}^{C}_x /\xi=j\}$$
will be called the \emph{level j of the simplicial complex} $\mathcal{K}^{C}_x$.


\begin{defn}
The \emph{interaction filtration} of $\mathcal{K}^{C}_x$ is the \emph{Hasse diagram} of the ordering generated out of every application of $\odot$ to the initial chain $x$.     
\end{defn}

\begin{example}
For the chain $x=G\otimes H\otimes K$ its interaction filtration contains 

$$G\otimes H\otimes K\leq G\odot H\otimes K\leq G\odot H\odot K$$ 
$$G\otimes H\otimes K\leq G\otimes H\odot K\leq G\odot H\odot K$$ 
$$G\otimes H\otimes K\leq G\odot K\otimes H\leq G\odot H\odot K$$ 
 
where: 
\begin{itemize}
    \item{$G\otimes H\otimes K \in \mathcal{K}_x[0]$}
    \item {$G\odot H\otimes K$, $G\odot K\otimes H$ and $G\otimes H\odot K \in \mathcal{K}_x[1]$}
    \item{$G\odot H\odot K \in \mathcal{K}_x[2]$.}
\end{itemize}
\end{example}
Observe that for a given $k$ the bottom of a poset filtration is a chain in the form $$G_{\pi(1)}\otimes\cdots\otimes G_{\pi(k)}$$ and not the empty set, this is one of the reasons why this sequence of complexes does not fit into the formal definition of filtration. On the other hand, the top chain $$G_{1}\odot\cdots\odot G_{k}$$ plays the role of the whole multicomplex (where every interaction has taken place).

\begin{lem} \label{conec}
For every $j \in [k]$:
\begin{itemize}
\item{the number of elements into $\mathcal{K}[j]$ is ${{k}\choose{j+1}}$} 
\item{every chain into $\mathcal{K}[j]$ consists of $k-j$ connected components.}      
\end{itemize}
   
\end{lem}

\subsection{Results on the interaction filtration}
A substantial difference with other filtrations such as those referred to in \cite{Akt}, is that in ours no more edges appear in every step of the complexification of the filtration. Moreover, to count the steps into the filtration as in \cite{Akt} we need something more than simply a chain of pictures (as in the totally ordered models). We will see this from the examples above, where the following Lemma is satisfied by any of them:
\begin{lem}
Given a simplicial complex $\mathcal{K}^{C}_x$ with $k=\mid C \mid$, its interaction filtration is p-folded, where $p$ equals to $$1+\stackrel[j=2]{k}{\sum}{{k}\choose{j}}=2^{k}-k$$
\end{lem}
To construct a full interaction filtration out of $k$ graphs we count the steps as done in the following example. Note that we are considering for simplicity the same colour for all the edges, while they will be coloured in next sections. Take $\delta \in \{0,...,p-1\}$. 

\begin{example}
For $k=3$ we have $p=4$, then for $\delta=0$, if $G \otimes H \otimes K$ has the form:
\vspace{1em}
\input{Firstfiltration.tikz}
\vspace{1em}

then for $\delta=1$ and $\delta=2$, $G \odot H \otimes K$ and $G \otimes H \odot K$ have respectively the forms:
\vspace{1em}
\input{Secondfiltration.tikz}
\vspace{1em}
\hspace{-1em}and finally, for $\delta=3$:
\vspace{1em}
\input{Thirdfiltration.tikz}
\end{example}
\subsection{Holes in a multicomplex}
As pointed in the introduction, the main target of this paper is to investigate which of the topological features are preserved through the process of merging graphs in $\mathcal{G}^{\otimes C}$. 

We are in particular interested in calculate the numbers of \emph{holes} appearing and remaining along the configurations of an interaction filtration. That is precisely the reason why one should care about filtrations: they provide a way to measure the lifespan of a hole of any dimension into a simplicial complex. This is the field of study of Persistent Homology as will be seen in the following section.

To introduce the concept of \emph{hole} in a simplicial complex one needs to define first that of \emph{chain group}, \emph{cycle} and \emph{boundary} of a simplicial complex (see \cite{What} for instance). This will be done in the next section where will be seen that, given a simplicial complex $\mathcal{K}$, the $0-holes$ are its connected components, the $1-holes$ are the cavities into $\mathcal{K}$ surrounded by edges, the $2-holes$ are voids enclosed by triangles, the $3-holes$ are surrounded by tetrahedrons etc. 

Note, however, that we exclude the cliques of any dimension: a hole is not a clique. That is: a $d-hole$, or a hole of dimension $d$, into a simplicial complex is a void which is not clique enclosed by a number of $d$-simplices.
\begin{example} \label{holes1}
Considering the following complexes
\vspace{1em}
\input{holes.tikz}
\vspace{1em}

\hspace{-1em}we conclude that: a) has one 0-hole and two 1-holes, b) has one 0-hole, one 1-hole and two 1-cliques while c) has one 0-hole, one 1-hole and one 3-clique. 
\end{example}
Observe that in our notation cliques are filled with a tone of grey while holes remain blank. Recall also that every time a clique is closed, it is not a hole any more. In the context of a filtration, as $\delta$ increases, the number and complexity of the holes also does:
\begin{example} \label{holes2}
For the filtration
\vspace{1em}
\input{Fourthfiltration.tikz}

\hspace{-1em}we can identify three 0-holes for $\delta=0$, two 0-holes and two 1-cliques for $\delta=1$, two 0-holes and one 1-clique for $\delta=2$ and one 0-hole and one 3-clique for $\delta=3$.

Notice that this corresponds to a poset with
$$G\otimes H\otimes K\leq G\otimes H\odot K \leq G\odot H\odot K$$ and $$G\otimes H\otimes K\leq G\odot H\otimes K \leq G\odot H\odot K$$
naming $G,H,K$ the graphs involved starting from above.
\end{example}

\section{Persistent Homology}
Persistent Homology serves to quantify some topological features along different configurations of the data, it appears as a way to calculate homology during a (discrete) scale of time. Our aim is to calculate it for our example in \cite{Sign} where every application of $\odot$ is seen as a step in a certain scale of time. 
\begin{defn}
The \emph{n-chain group} of a simplicial complex $\mathcal{K}$, denoted by $C_n(\mathcal{K},\mathbb{F})$, is the vector space generated by the n-simplices of $\mathcal{K}$ with coefficients in a field $\mathbb{F}$.
\end{defn}
That is, the elements of $C_n(\mathcal{K},\mathbb{F})$ are sums of simplices in $\mathcal{K}$. For Persistent Homology it is common to use $\mathbb{F}=\mathbb{Z}_2$.
\begin{defn}
The \emph{boundary operators} of $C_n(\mathcal{K},\mathbb{F})$ are the linear maps $$\partial:C_n(\mathcal{K},\mathbb{F}) \rightarrow C_{n-1}(\mathcal{K},\mathbb{F})$$ such that $$\partial(v_0,...,v_n)=\stackrel[j=0]{n}{\sum}(-1)^{j}(v_{0},...,v_{\{j-1\}},v_{\{j+1\}},...,v_{n})$$ where we denote by $(v_0,...,v_n)$ the (n+1)-simplice in $\mathcal{K}$ composed by those vertices. For that (n+1)-simplice, every tuple $(v_{0},...,v_{\{j-1\}},v_{\{j+1\}},...,v_{n})$ is a n-face.
\end{defn}
\subsection{Betti numbers}
The number of n-holes in a configuration of a filtration is called its \emph{nth Betti number}. Formally (and see \cite{What} for an extensive explanation):

\begin{defn}
    The \emph{n-dimensional Betti number}, denoted by $\beta_n$, is the dimension of the quotient $${Ker \partial_n} \over {Im \partial_{n+1}}$$
\end{defn}

\begin{example} For the complexes in Example \ref{holes1} we have $\beta_0=1$, $\beta_1=2$ in a), $\beta_0=1$, $\beta_1=1$ in b) and $\beta_0=1$, $\beta_1=1$ in c).

For the complexes in Example \ref{holes2} we have $\beta_0=3$ when $\delta=0$, $\beta_0=2$ when $\delta=1$, $\beta_0=2$ when $\delta=2$ and $\beta_0=1$ when $\delta=3$. 

From Example \ref{holes2} we see that, when the number of $\delta$ increases, the Betti number of a given dimension decreases, this is coherent with the algorithm provided in \cite{Incremental}. Something more is happening, however, as be seen in the following.
   
\end{example}
\subsection{The incremental method} 
Observe that Lemma \ref{conec} says that every chain into a certain $\mathcal{K}^{C}[j]$ has the same $\beta_0$. It is a natural question asking how to calculate the first three dimensional Betti numbers of every chain into the same $\mathcal{K}^{C}[j]$, this can be done by modifying the \emph{incremental algorithm} designed in \cite{Incremental}.

The results given in \cite{Incremental} are different from ours since two crucial aspects of the filtration considered are also different. On one hand, for us simplices appear in some steps of the filtration only when they correspond to a clique in the underlying graph. Secondly, we are working here multigraphs rather than graphs. These two facts make the formula to calculate the Betti numbers more complicated as will be seen. It appears, however, as a sort of generalization of the one given in \cite{Incremental}. 

To obtain the general formula, that working for multicomplexes, we first need some preliminary results. The first lemma corresponds to the \emph{incremental method} introduced in \cite{Incremental} for \emph{simple complexes} (that is, those defined out of simple graphs) and built using a \emph{usual filtration}. 

\begin{lem}\label{lema1}
Let $\mathcal{K}$ and $\mathcal{K'}$ simplicial complexes such that $\mathcal{K'}$ is obtained by adding the d-simplex $\sigma$ to $\mathcal{K}$ in a filtration process. Then:
\begin{itemize}
\item{if $\sigma$ belongs to a d-cycle in $\mathcal{K'}$: $$\beta_d(\mathcal{K'})=\beta_{d}(\mathcal{K})+1$$}
\item{if $\sigma$ does not belong to a d-cycle in $\mathcal{K'}$: $$\beta_{d-1}(\mathcal{K'})=\beta_{d-1}(\mathcal{K})-1$$}
\end{itemize}
\end{lem}
Let us finish these preliminar results by observing that \emph{duplications generate holes}:
\begin{lem}\label{lema2}
Let $\mathcal{K}$ and $\mathcal{K'}$ multicomplexes such that $\mathcal{K} \subseteq \mathcal{K'}$. Each d-dimensional duplication involved in a d-cycle of $\mathcal{K}$ generates a (d+1)-hole in $\mathcal{K'}$.  
\end{lem}

Now we are in position of introducing analogous formulas to those of \cite{Incremental} for multi- complexes built using the interaction filtration. The formulae in the following result appear as a sort of \emph{Inclusion-Exclusion Principle} for multicomplexes.

\begin{prop}\label{formulae}

    Let be $G,H$ multigraphs and $\mathcal{K}_{G\odot H}$ its associated multicomplex. Considering its interaction filtration, the Betti numbers of the configuration $\mathcal{K}_{G\odot H}$ are given by the following formulae:
    $$\beta_{1}(\mathcal{K}_{G\odot H})=max(\beta_{1}(\mathcal{K}_{G}),\beta_{1}(\mathcal{K}_{H}))+max(n_{G},n_{H})-min(p_G,p_H)-cl$$
    $$\beta_{2}(\mathcal{K}_{G\odot H})=max(\beta_{2}(\mathcal{K}_{G}),\beta_{2}(\mathcal{K}_{H}))+max(n_{G},n_{H})-min(p_G,p_H)-cl+dup$$
    where each $n$ is the number of simplices added that close a d-cycle, p is the number of simplices added that do not close a d-cycle, cl is the number of cliques created in each case and dup the number of duplications for $d=1$.

\end{prop}

\begin{proof}
The idea is to follow the algorithm designed in \cite{Incremental} by adapting it to the particularities of merging multigraphs. That is, when adding a new simplice $\sigma$ in an interaction filtration through $\odot$ we must distinguish among three possibilities: $\sigma$ closes a \emph{d-cycle} in the new configuration, $\sigma$ does not close a \emph{d-cycle} in the new configuration and $\sigma$ duplicates an existing simplice in the previous configuration.

It is crucial to observe that the interaction filtration is a usual filtration (as that of \cite{Incremental} for instance) for which we add several simplices in a single step. This allows us to apply the formulas in \ref{lema1} to our context. We should be careful, however, with the fact that, from Definition \ref{cliquesimplex}, we are closing every clique: since we are counting holes, we have to remove every clique appearing after a new inclusion of a simplice. This points to the \emph{Inclusion-Exclusion} looking of our formulas. 

Therefore, in the first case, and adding one by one the new simplices that close a \emph{d-cycle}, we must sum the number of them to the Betti number of the previous configuration. Now is time to substract the cycles that convert into cliques, since we are counting \emph{d-holes}.\footnote{In the terms $min(p_G,p_H)$ of the formulas there is a flavour to the \emph{Morgan Laws} in a Logic system.} This is a substantial difference with the calculations made in \cite{Incremental}, where cliques are not immediately simplices as in our case and where only one simplice is added in every configuration of the filtration. This case includes parallel simplices, that is, simplices with the same vertices in the underlying multigraph. 

For the simplices that do not close a \emph{d-cycle}, just do the same process: step by step add those simplices and calculate accordingly. Observe that the term $dup$ appears as a consequence of Lemma \ref{lema2}.
\end{proof}

\begin{example}
A) For the composition
    \input{betti1.tikz}
    
\hspace{-1em}we calculate:
\begin{itemize}
    \item{for $d=1$ we have $\beta_{1}(\mathcal{K}_{G})=1$ and $\beta_{1}(\mathcal{K}_{H})=0$ and the assignations: $$n_G=0,cl=0,p_G=1$$ as well as $$n_H=0,cl=0,p_H=1$$ and then $\beta_{1}(\mathcal{K}_{G\odot H})=1$}
    \item{for $d=2$ we have $\beta_{2}(\mathcal{K}_{G})=0$ and $\beta_{2}(\mathcal{K}_{H})=0$ and the assignations: $$n_G=0,cl=0,p_G=2$$ as well as $$n_H=0,cl=0,p_H=1$$ and since $dup=1$ we obtain $\beta_{2}(\mathcal{K}_{G\odot H})=0$.}
\end{itemize}

B) For the composition
    \input{betti2.tikz}
    
we calculate:
\begin{itemize}
    \item{for $d=1$ we have $\beta_{1}(\mathcal{K}_{G})=1$ and $\beta_{1}(\mathcal{K}_{H})=0$ and the assignations: $$n_G=2,cl=1,p_G=0$$ as well as $$n_H=1,cl=1,p_H=1$$ and then $\beta_{1}(\mathcal{K}_{G\odot H})=1$}
    \item{for $d=2$ we have $\beta_{2}(\mathcal{K}_{G})=0$ and $\beta_{2}(\mathcal{K}_{H})=0$ and the assignations: $$n_G=0,cl=0,p_G=0$$ as well as $$n_H=0,cl=0,p_H=0$$ and since $dup=1$ we obtain $\beta_{2}(\mathcal{K}_{G\odot H})=1$.}
\end{itemize}

C) For the composition (where we have not filled the tetrahedron)
    \input{betti3.tikz}
    
we calculate:
\begin{itemize}
    \item{for $d=1$ we have $\beta_{1}(\mathcal{K}_{G})=1$ and $\beta_{1}(\mathcal{K}_{H})=0$ and the assignations: $$n_G=2,cl=2,p_G=0$$ as well as $$n_H=0,cl=2,p_H=0$$ and then $\beta_{1}(\mathcal{K}_{G\odot H})=1$}
    \item{for $d=2$ we have $\beta_{2}(\mathcal{K}_{G})=0$ and $\beta_{2}(\mathcal{K}_{H})=0$ and the assignations: $$n_G=0,cl=2,p_G=0$$ as well as $$n_H=2,cl=2,p_H=0$$ and since $dup=1$ we obtain $\beta_{2}(\mathcal{K}_{G\odot H})=1$.}
\end{itemize}
\end{example}

\subsection*{Acknowledgements} The author wants to gratefully acknowledge the help of Prof. Orlando Galdames Bravo, whose valuable suggestions made the second section of this paper much deeper and clearer.

\end{document}

%% file: exampletensor.tikz
\tikzset{every picture/.style={line width=0.75pt}} 
\vspace{1em}
\hspace{3em}
\begin{tikzpicture}[x=0.75pt,y=0.75pt,yscale=-1,xscale=1]

\draw [color={rgb, 255:red, 74; green, 144; blue, 226 }  ,draw opacity=1 ][line width=2.25]    (161.31,214.57) .. controls (158.67,200.78) and (170.24,174.27) .. (197.52,167.49) ;
\draw [color={rgb, 255:red, 74; green, 144; blue, 226 }  ,draw opacity=1 ][line width=2.25]    (165.94,215.63) -- (200.99,171.73) ;
\draw [color={rgb, 255:red, 208; green, 2; blue, 27 }  ,draw opacity=1 ][line width=2.25]    (167.93,219.87) .. controls (190.58,220.08) and (207.94,182.76) .. (205.95,171.73) ;
\draw [color={rgb, 255:red, 74; green, 144; blue, 226 }  ,draw opacity=1 ][line width=2.25]    (208.1,166.85) -- (225.79,166.94) -- (252.08,167.06) ;
\draw [color={rgb, 255:red, 126; green, 211; blue, 33 }  ,draw opacity=1 ][line width=2.25]    (291.42,221.99) -- (316.55,172.58) ;
\draw [color={rgb, 255:red, 248; green, 231; blue, 28 }  ,draw opacity=1 ][line width=2.25]    (295.89,223.26) -- (348.3,174.27) ;
\draw [color={rgb, 255:red, 126; green, 211; blue, 33 }  ,draw opacity=1 ][line width=2.25]    (376.07,218.81) -- (353.25,173.43) ;
\draw [color={rgb, 255:red, 248; green, 231; blue, 28 }  ,draw opacity=1 ][line width=2.25]    (357.39,170.24) .. controls (376.07,166.64) and (399.38,207.36) .. (382.02,217.96) ;
\draw [color={rgb, 255:red, 74; green, 144; blue, 226 }  ,draw opacity=1 ][line width=2.25]    (417.11,211.9) .. controls (414.47,198.12) and (426.04,171.61) .. (453.32,164.82) ;
\draw [color={rgb, 255:red, 74; green, 144; blue, 226 }  ,draw opacity=1 ][line width=2.25]    (421.74,212.96) -- (456.79,169.06) ;
\draw [color={rgb, 255:red, 208; green, 2; blue, 27 }  ,draw opacity=1 ][line width=2.25]    (423.73,217.2) .. controls (446.38,217.42) and (463.73,180.09) .. (461.75,169.06) ;
\draw [color={rgb, 255:red, 74; green, 144; blue, 226 }  ,draw opacity=1 ][line width=2.25]    (463.9,164.18) -- (481.59,164.27) -- (507.88,164.39) ;
\draw [color={rgb, 255:red, 126; green, 211; blue, 33 }  ,draw opacity=1 ][line width=2.25]    (525.97,216.33) -- (551.1,166.91) ;
\draw [color={rgb, 255:red, 248; green, 231; blue, 28 }  ,draw opacity=1 ][line width=2.25]    (530.43,217.6) -- (582.84,168.61) ;
\draw [color={rgb, 255:red, 126; green, 211; blue, 33 }  ,draw opacity=1 ][line width=2.25]    (610.61,213.15) -- (587.8,167.76) ;
\draw [color={rgb, 255:red, 248; green, 231; blue, 28 }  ,draw opacity=1 ][line width=2.25]    (591.93,164.58) .. controls (610.61,160.97) and (633.92,201.69) .. (616.56,212.3) ;

\draw (157.97,218.54) node [anchor=north west][inner sep=0.75pt]  [font=\Large,xscale=0.5,yscale=0.5]  {$1$};
\draw (198.13,158.4) node [anchor=north west][inner sep=0.75pt]  [font=\Large,xscale=0.5,yscale=0.5]  {$2$};
\draw (255.2,160.83) node [anchor=north west][inner sep=0.75pt]  [font=\Large,xscale=0.5,yscale=0.5]  {$3$};
\draw (260.2,182.88) node [anchor=north west][inner sep=0.75pt]  [font=\Huge,xscale=0.5,yscale=0.5]  {$\otimes $};
\draw (286.42,220.91) node [anchor=north west][inner sep=0.75pt]  [font=\Large,xscale=0.5,yscale=0.5]  {$1$};
\draw (311.71,159.25) node [anchor=north west][inner sep=0.75pt]  [font=\Large,xscale=0.5,yscale=0.5]  {$2$};
\draw (345.93,159.68) node [anchor=north west][inner sep=0.75pt]  [font=\Large,xscale=0.5,yscale=0.5]  {$3$};
\draw (373.04,220.88) node [anchor=north west][inner sep=0.75pt]  [font=\Large,xscale=0.5,yscale=0.5]  {$4$};
\draw (394.41,187.88) node [anchor=north west][inner sep=0.75pt]  [font=\Huge,xscale=0.5,yscale=0.5]  {$=$};
\draw (413.77,216.88) node [anchor=north west][inner sep=0.75pt]  [font=\Large,xscale=0.5,yscale=0.5]  {$1$};
\draw (453.93,155.73) node [anchor=north west][inner sep=0.75pt]  [font=\Large,xscale=0.5,yscale=0.5]  {$2$};
\draw (510,156.16) node [anchor=north west][inner sep=0.75pt]  [font=\Large,xscale=0.5,yscale=0.5]  {$3$};
\draw (521.96,217.24) node [anchor=north west][inner sep=0.75pt]  [font=\Large,xscale=0.5,yscale=0.5]  {$1$};
\draw (548.25,154.59) node [anchor=north west][inner sep=0.75pt]  [font=\Large,xscale=0.5,yscale=0.5]  {$2$};
\draw (582.48,155.01) node [anchor=north west][inner sep=0.75pt]  [font=\Large,xscale=0.5,yscale=0.5]  {$3$};
\draw (612.61,216.55) node [anchor=north west][inner sep=0.75pt]  [font=\Large,xscale=0.5,yscale=0.5]  {$4$};

\end{tikzpicture}

%% file: odot.tikz
\tikzset{every picture/.style={line width=0.75pt}} 
\hspace{3em}
\begin{tikzpicture}[x=0.75pt,y=0.75pt,yscale=-1,xscale=1]

\draw [color={rgb, 255:red, 74; green, 144; blue, 226 }  ,draw opacity=1 ][line width=2.25]    (24,143) .. controls (-2,130.42) and (35,75.42) .. (49,82.42) ;
\draw [color={rgb, 255:red, 208; green, 2; blue, 27 }  ,draw opacity=1 ][line width=2.25]    (63,86.42) .. controls (84,98.42) and (57,142.42) .. (29,143.42) ;
\draw [color={rgb, 255:red, 74; green, 144; blue, 226 }  ,draw opacity=1 ][line width=2.25]    (53,88.42) -- (27,138.42) ;
\draw [color={rgb, 255:red, 74; green, 144; blue, 226 }  ,draw opacity=1 ][line width=2.25]    (63,80.42) -- (112,80.42) ;
\draw [color={rgb, 255:red, 126; green, 211; blue, 33 }  ,draw opacity=1 ][line width=2.25]    (188,87.42) -- (163,142.42) ;
\draw [color={rgb, 255:red, 126; green, 211; blue, 33 }  ,draw opacity=1 ][line width=2.25]    (238,89) -- (268,139.42) ;
\draw [color={rgb, 255:red, 248; green, 231; blue, 28 }  ,draw opacity=1 ][line width=2.25]    (164,149.42) -- (224,86.42) ;
\draw [color={rgb, 255:red, 248; green, 231; blue, 28 }  ,draw opacity=1 ][line width=2.25]    (241,78) .. controls (287,86.42) and (277,108.42) .. (275,140.42) ;
\draw [color={rgb, 255:red, 126; green, 211; blue, 33 }  ,draw opacity=1 ][line width=2.25]    (360,87.42) -- (335,142.42) ;
\draw [color={rgb, 255:red, 74; green, 144; blue, 226 }  ,draw opacity=1 ][line width=2.25]    (370,79.42) -- (406,79.42) ;
\draw [color={rgb, 255:red, 126; green, 211; blue, 33 }  ,draw opacity=1 ][line width=2.25]    (423,90) -- (453,140.42) ;
\draw [color={rgb, 255:red, 248; green, 231; blue, 28 }  ,draw opacity=1 ][line width=2.25]    (339,152) -- (409,87.42) ;
\draw [color={rgb, 255:red, 248; green, 231; blue, 28 }  ,draw opacity=1 ][line width=2.25]    (426,79) .. controls (472,87.42) and (462,109.42) .. (460,141.42) ;
\draw [color={rgb, 255:red, 74; green, 144; blue, 226 }  ,draw opacity=1 ][line width=2.25]    (325,147.5) .. controls (299,134.92) and (340,76.42) .. (354,83.42) ;
\draw [color={rgb, 255:red, 74; green, 144; blue, 226 }  ,draw opacity=1 ][line width=2.25]    (329,144) .. controls (319,123.5) and (348,90.58) .. (355,88.5) ;
\draw [color={rgb, 255:red, 208; green, 2; blue, 27 }  ,draw opacity=1 ][line width=2.25]    (366,90.5) .. controls (364,105.5) and (364.88,115.41) .. (359,125.5) .. controls (353.12,135.59) and (345,141.5) .. (341,146.42) ;

\draw (17,146) node [anchor=north west][inner sep=0.75pt]   [align=left] {1};
\draw (51,71) node [anchor=north west][inner sep=0.75pt]   [align=left] {2};
\draw (117,71.4) node [anchor=north west][inner sep=0.75pt]    {$3$};
\draw (132,97.4) node [anchor=north west][inner sep=0.75pt]  [font=\large]  {$\odot $};
\draw (287,104.4) node [anchor=north west][inner sep=0.75pt]  [font=\large]  {$=$};
\draw (152,145) node [anchor=north west][inner sep=0.75pt]   [align=left] {1};
\draw (186,74) node [anchor=north west][inner sep=0.75pt]   [align=left] {2};
\draw (227,74.4) node [anchor=north west][inner sep=0.75pt]    {$3$};
\draw (265,141.4) node [anchor=north west][inner sep=0.75pt]    {$4$};
\draw (324,145) node [anchor=north west][inner sep=0.75pt]   [align=left] {1};
\draw (355,73) node [anchor=north west][inner sep=0.75pt]   [align=left] {2};
\draw (411,75.4) node [anchor=north west][inner sep=0.75pt]    {$3$};
\draw (450,142.4) node [anchor=north west][inner sep=0.75pt]    {$4$};

\end{tikzpicture}

%% file: Firstclique.tikz
\tikzset{every picture/.style={line width=0.6pt}} 
\hspace{1em}

\begin{tikzpicture}[x=0.65pt,y=0.65pt,yscale=-1,xscale=1]

\draw [line width=3]    (27.09,107.99) -- (96.5,108) ;
\draw [line width=3]    (99.51,178) -- (100.52,112) ;
\draw [line width=3]    (94.48,179) -- (25.08,114.99) ;
\draw [line width=3]    (102.54,94) .. controls (99.52,77) and (114.62,43) .. (167.93,60.01) ;
\draw [line width=3]    (209.49,100.43) -- (278.48,99.53) ;
\draw [line width=3]    (283.45,173.47) -- (282.53,103.47) ;
\draw [line width=3]    (278.47,175.53) -- (207.58,107.45) ;
\draw [line width=3]    (276.57,106.55) -- (248.95,135.92) ;
\draw [line width=3]    (241.06,144.02) -- (207.48,176.46) ;
\draw [line width=3]    (201.43,172.54) -- (201.62,110.53) ;
\draw [line width=3]    (209.56,182.43) -- (278.55,181.53) ;
\draw [line width=3]    (429.49,93) .. controls (426.48,76) and (441.57,42) .. (494.88,59.01) ;
\draw [line width=3]    (536.44,99.43) -- (605.44,98.53) ;
\draw [line width=3]    (610.4,172.47) -- (609.49,102.47) ;
\draw [line width=3]    (605.43,174.53) -- (534.53,106.45) ;
\draw [line width=3]    (604.53,105.55) -- (574,137) ;
\draw [line width=3]    (568.02,143.02) -- (534.44,175.46) ;
\draw [line width=3]    (528.39,171.54) -- (528.57,109.53) ;
\draw [line width=3]    (536.52,181.43) -- (605.51,180.53) ;
\draw [line width=3]    (356.05,103.99) -- (425.46,104) ;
\draw [line width=3]    (428.46,174) -- (429.48,108) ;
\draw [line width=3]    (423.43,177) -- (354.03,112.99) ;
\draw  [fill={rgb, 255:red, 155; green, 155; blue, 155 }  ,fill opacity=0.34 ][dash pattern={on 0.84pt off 2.51pt}][line width=0.75]  (358,109) -- (423,170) -- (423.73,109.79) -- cycle ;
\draw  [fill={rgb, 255:red, 155; green, 155; blue, 155 }  ,fill opacity=1 ][dash pattern={on 0.84pt off 2.51pt}][line width=0.75]  (606.68,167.42) -- (606.94,109.02) -- (575.79,138.67) -- cycle ;
\draw  [fill={rgb, 255:red, 155; green, 155; blue, 155 }  ,fill opacity=1 ][dash pattern={on 0.84pt off 2.51pt}][line width=0.75]  (601,102) -- (536.29,102.86) -- (569.07,134.07) -- cycle ;
\draw  [fill={rgb, 255:red, 155; green, 155; blue, 155 }  ,fill opacity=1 ][dash pattern={on 0.84pt off 2.51pt}][line width=0.75]  (532.23,109.55) -- (531.86,172.91) -- (563,140) -- cycle ;
\draw  [fill={rgb, 255:red, 155; green, 155; blue, 155 }  ,fill opacity=1 ][dash pattern={on 0.84pt off 2.51pt}][line width=0.75]  (540.06,177.3) -- (600,177) -- (569.15,147.85) -- cycle ;

\draw (13.04,97.99) node [anchor=north west][inner sep=0.75pt]  [rotate=-0.01,xscale=0.75,yscale=0.75] [align=left] {1};
\draw (96.52,182) node [anchor=north west][inner sep=0.75pt]  [rotate=-0.01,xscale=0.75,yscale=0.75] [align=left] {2};
\draw (96.53,96) node [anchor=north west][inner sep=0.75pt]  [rotate=-0.01,xscale=0.75,yscale=0.75] [align=left] {3};
\draw (169.97,54.01) node [anchor=north west][inner sep=0.75pt]  [rotate=-0.01,xscale=0.75,yscale=0.75] [align=left] {4};
\draw (195.36,90.61) node [anchor=north west][inner sep=0.75pt]  [rotate=-359.25,xscale=0.75,yscale=0.75] [align=left] {5};
\draw (279.45,173.52) node [anchor=north west][inner sep=0.75pt]  [rotate=-359.25,xscale=0.75,yscale=0.75] [align=left] {8};
\draw (277.32,87.54) node [anchor=north west][inner sep=0.75pt]  [rotate=-359.25,xscale=0.75,yscale=0.75] [align=left] {7};
\draw (195.42,171.62) node [anchor=north west][inner sep=0.75pt]  [rotate=-359.25,xscale=0.75,yscale=0.75] [align=left] {6};
\draw (496.93,53.01) node [anchor=north west][inner sep=0.75pt]  [rotate=-0.01,xscale=0.75,yscale=0.75] [align=left] {4};
\draw (522.31,89.61) node [anchor=north west][inner sep=0.75pt]  [rotate=-359.25,xscale=0.75,yscale=0.75] [align=left] {5};
\draw (606.4,172.52) node [anchor=north west][inner sep=0.75pt]  [rotate=-359.25,xscale=0.75,yscale=0.75] [align=left] {8};
\draw (604.28,86.54) node [anchor=north west][inner sep=0.75pt]  [rotate=-359.25,xscale=0.75,yscale=0.75] [align=left] {7};
\draw (522.37,170.62) node [anchor=north west][inner sep=0.75pt]  [rotate=-359.25,xscale=0.75,yscale=0.75] [align=left] {6};
\draw (342,93.99) node [anchor=north west][inner sep=0.75pt]  [rotate=-0.01,xscale=0.75,yscale=0.75] [align=left] {1};
\draw (425.48,178) node [anchor=north west][inner sep=0.75pt]  [rotate=-0.01,xscale=0.75,yscale=0.75] [align=left] {2};
\draw (425.49,92) node [anchor=north west][inner sep=0.75pt]  [rotate=-0.01,xscale=0.75,yscale=0.75] [align=left] {3};

\end{tikzpicture}

%% file: Secondclique.tikz
\tikzset{every picture/.style={line width=0.75pt}} 
\hspace{6em}
\begin{tikzpicture}[x=0.75pt,y=0.75pt,yscale=-1,xscale=1]

\draw [line width=3]    (259.54,100) .. controls (256.52,83) and (271.62,49) .. (324.93,66.01) ;
\draw [line width=3]    (366.49,106.43) -- (435.48,105.53) ;
\draw [line width=3]    (440.45,179.47) -- (439.53,109.47) ;
\draw [line width=3]    (358.43,178.54) -- (358.62,116.53) ;
\draw [line width=3]    (366.56,188.43) -- (435.55,187.53) ;
\draw [line width=3]    (185.05,109.99) -- (254.46,110) ;
\draw [line width=3]    (257.46,180) -- (258.48,114) ;
\draw [line width=3]    (252.43,183) -- (183.03,118.99) ;
\draw  [fill={rgb, 255:red, 155; green, 155; blue, 155 }  ,fill opacity=1 ][dash pattern={on 0.84pt off 2.51pt}][line width=0.75]  (187,115) -- (252,176) -- (252.73,115.79) -- cycle ;

\draw (326.97,60.01) node [anchor=north west][inner sep=0.75pt]  [rotate=-0.01,xscale=0.75,yscale=0.75] [align=left] {4};
\draw (352.36,96.61) node [anchor=north west][inner sep=0.75pt]  [rotate=-359.25,xscale=0.75,yscale=0.75] [align=left] {5};
\draw (436.45,179.52) node [anchor=north west][inner sep=0.75pt]  [rotate=-359.25,xscale=0.75,yscale=0.75] [align=left] {8};
\draw (434.32,93.54) node [anchor=north west][inner sep=0.75pt]  [rotate=-359.25,xscale=0.75,yscale=0.75] [align=left] {7};
\draw (352.42,177.62) node [anchor=north west][inner sep=0.75pt]  [rotate=-359.25,xscale=0.75,yscale=0.75] [align=left] {6};
\draw (171,99.99) node [anchor=north west][inner sep=0.75pt]  [rotate=-0.01,xscale=0.75,yscale=0.75] [align=left] {1};
\draw (254.48,184) node [anchor=north west][inner sep=0.75pt]  [rotate=-0.01,xscale=0.75,yscale=0.75] [align=left] {2};
\draw (254.49,98) node [anchor=north west][inner sep=0.75pt]  [rotate=-0.01,xscale=0.75,yscale=0.75] [align=left] {3};

\end{tikzpicture}

%% file: Thirdclique.tikz
\tikzset{every picture/.style={line width=0.75pt}} 
\hspace{4em}
\begin{tikzpicture}[x=0.65pt,y=0.65pt,yscale=-1,xscale=1]

\draw [line width=3]    (90,51) -- (180,51) ;
\draw [line width=3]    (251,52) -- (341,52) ;
\draw [line width=3]    (344,57) -- (343,138) ;
\draw [line width=3]    (425,54) -- (515,54) ;
\draw [line width=3]    (423,59) -- (513,142) ;
\draw [line width=3]    (519,58) -- (518,139) ;
\draw [line width=3]    (90,58) -- (180,141) ;
\draw [line width=3]    (422,49) .. controls (445,2) and (503,16) .. (517,48) ;
\draw  [draw opacity=0][fill={rgb, 255:red, 155; green, 155; blue, 155 }  ,fill opacity=0.52 ][dash pattern={on 0.84pt off 2.51pt}] (426.78,51.17) .. controls (426.78,51.04) and (426.78,50.91) .. (426.78,50.79) .. controls (427,35.2) and (446.55,22.98) .. (470.43,23.5) .. controls (493.03,23.99) and (511.42,35.7) .. (513.15,50.16) -- (470.03,51.72) -- cycle ; \draw  [dash pattern={on 0.84pt off 2.51pt}] (426.78,51.17) .. controls (426.78,51.04) and (426.78,50.91) .. (426.78,50.79) .. controls (427,35.2) and (446.55,22.98) .. (470.43,23.5) .. controls (493.03,23.99) and (511.42,35.7) .. (513.15,50.16) ;  
\draw  [fill={rgb, 255:red, 128; green, 128; blue, 128 }  ,fill opacity=1 ][dash pattern={on 0.84pt off 2.51pt}][line width=0.75]  (429,58) -- (513,136) -- (513.93,59.02) -- cycle ;

\draw (198,78.4) node [anchor=north west][inner sep=0.75pt]  [font=\Large,xscale=0.75,yscale=0.75]  {$\bigodot $};
\draw (380,94.4) node [anchor=north west][inner sep=0.75pt]  [font=\fontsize{4.41em}{5.29em}\selectfont,xscale=0.75,yscale=0.75]  {$=$};
\draw (79,42.4) node [anchor=north west][inner sep=0.75pt]  [xscale=0.75,yscale=0.75]  {$1$};
\draw (182,134.4) node [anchor=north west][inner sep=0.75pt]  [xscale=0.75,yscale=0.75]  {$2$};
\draw (181,43.4) node [anchor=north west][inner sep=0.75pt]  [xscale=0.75,yscale=0.75]  {$3$};
\draw (242,43.4) node [anchor=north west][inner sep=0.75pt]  [xscale=0.75,yscale=0.75]  {$1$};
\draw (337,136.4) node [anchor=north west][inner sep=0.75pt]  [xscale=0.75,yscale=0.75]  {$2$};
\draw (342,41.4) node [anchor=north west][inner sep=0.75pt]  [xscale=0.75,yscale=0.75]  {$3$};
\draw (519,42.4) node [anchor=north west][inner sep=0.75pt]  [xscale=0.75,yscale=0.75]  {$3$};
\draw (410,44.4) node [anchor=north west][inner sep=0.75pt]  [xscale=0.75,yscale=0.75]  {$1$};
\draw (514,138.4) node [anchor=north west][inner sep=0.75pt]  [xscale=0.75,yscale=0.75]  {$2$};

\end{tikzpicture}

%% file: Forthclique.tikz
\tikzset{every picture/.style={line width=0.75pt}} 
\hspace{4em}
\begin{tikzpicture}[x=0.6pt,y=0.6pt,yscale=-1,xscale=1]

\draw [line width=3]    (68,52) -- (125,93) ;
\draw [line width=3]    (130,108) -- (129,183) ;
\draw [line width=3]    (134,93) -- (187,58) ;
\draw [line width=3]    (249,56) -- (366,56) ;
\draw [line width=3]    (245,65) -- (310,187) ;
\draw [line width=3]    (315,97) -- (368,62) ;
\draw [line width=3]    (372,67) -- (316,188) ;
\draw [line width=3]    (456,68) -- (516,191) ;
\draw [line width=3]    (617,66) -- (528,189) ;
\draw [line width=3]    (523,113) -- (522,188) ;
\draw [line width=3]    (460,61) -- (517,102) ;
\draw [line width=3]    (529,101) -- (606,63) ;
\draw [line width=3]    (463,54) .. controls (494,11) and (567,7) .. (610,56) ;
\draw  [fill={rgb, 255:red, 155; green, 155; blue, 155 }  ,fill opacity=1 ][dash pattern={on 0.84pt off 2.51pt}] (462.17,68.84) -- (517,183) -- (518.74,110.89) -- cycle ;
\draw  [fill={rgb, 255:red, 155; green, 155; blue, 155 }  ,fill opacity=1 ][dash pattern={on 0.84pt off 2.51pt}] (612.15,61.94) -- (525.17,184.45) -- (531,106) -- cycle ;
\draw  [draw opacity=0][fill={rgb, 255:red, 74; green, 74; blue, 74 }  ,fill opacity=1 ][dash pattern={on 0.84pt off 2.51pt}] (528.04,99.23) .. controls (527.91,99.02) and (527.78,98.8) .. (527.67,98.58) .. controls (521.99,87.89) and (533.34,71.06) .. (553.01,60.98) .. controls (572.69,50.9) and (593.23,51.4) .. (598.9,62.09) .. controls (599.11,62.48) and (599.29,62.87) .. (599.45,63.28) -- (563.28,80.34) -- cycle ; \draw  [dash pattern={on 0.84pt off 2.51pt}] (528.04,99.23) .. controls (527.91,99.02) and (527.78,98.8) .. (527.67,98.58) .. controls (521.99,87.89) and (533.34,71.06) .. (553.01,60.98) .. controls (572.69,50.9) and (593.23,51.4) .. (598.9,62.09) .. controls (599.11,62.48) and (599.29,62.87) .. (599.45,63.28) ;  
\draw  [color={rgb, 255:red, 155; green, 155; blue, 155 }  ,draw opacity=1 ][line width=5.25] [line join = round][line cap = round] (489,62) .. controls (494.89,56.11) and (501.77,43.37) .. (511,43) .. controls (518.99,42.68) and (527.01,42.6) .. (535,43) .. controls (536.37,43.07) and (532.3,43.57) .. (531,44) .. controls (528.3,44.9) and (525.47,45.59) .. (523,47) .. controls (520.75,48.29) and (509.05,58.1) .. (511,62) .. controls (511.47,62.94) and (513.2,61.68) .. (514,61) .. controls (522.03,54.2) and (529.6,47.3) .. (538,41) .. controls (539.44,39.92) and (542.49,38.76) .. (544,38) .. controls (545.23,37.39) and (548.91,35.97) .. (548,37) .. controls (539.45,46.77) and (523.38,57.09) .. (521,69) .. controls (519.61,75.94) and (517.58,81.84) .. (515,87) .. controls (514.2,88.61) and (514.61,92.8) .. (513,92) .. controls (511.04,91.02) and (509.98,84) .. (508,84) .. controls (499.3,84) and (493.84,71.56) .. (487,67) .. controls (482.48,63.99) and (477.97,62.38) .. (474,60) .. controls (472.33,59) and (467.83,58.55) .. (469,57) .. controls (475.13,48.83) and (495.79,34.07) .. (505,31) .. controls (509.53,29.49) and (516.81,30.32) .. (521,30) .. controls (533.62,29.03) and (544.38,25.69) .. (555,31) .. controls (562.9,34.95) and (571.56,36.28) .. (579,40) .. controls (579.52,40.26) and (588.13,43.93) .. (588,44) .. controls (581.93,47.04) and (566.99,46.27) .. (561,46) .. controls (560.04,45.96) and (554.18,43.9) .. (554,44) .. controls (551.51,45.42) and (549.39,47.41) .. (547,49) .. controls (538.68,54.55) and (532.44,60.85) .. (526,66) .. controls (523.97,67.63) and (522.33,69.84) .. (520,71) .. controls (519.24,71.38) and (522.93,78.29) .. (518,75) .. controls (515.76,73.51) and (515.72,70.06) .. (514,68) .. controls (510.68,64.02) and (507.72,59.15) .. (503,57) .. controls (499.77,55.53) and (506,69.55) .. (506,69) .. controls (506,65.61) and (504.77,58.47) .. (501,58) .. controls (498.68,57.71) and (496.28,57.49) .. (494,58) .. controls (493.08,58.2) and (491.33,59.33) .. (492,60) .. controls (495.6,63.6) and (503.71,54.41) .. (502,51) .. controls (501.58,50.16) and (500.23,52.09) .. (500,53) .. controls (499.51,54.94) and (499.51,57.06) .. (500,59) .. controls (500.71,61.86) and (511.25,63.99) .. (512,61) .. controls (512.4,59.4) and (513.69,51.44) .. (509,53) .. controls (505.4,54.2) and (510.38,61.86) .. (513,54) .. controls (514.72,48.83) and (512.26,45) .. (507,45) .. controls (504.64,45) and (505.67,49.67) .. (506,52) .. controls (506.4,54.81) and (509.3,50.66) .. (510,53) .. controls (511.17,56.88) and (507.94,65) .. (512,65) .. controls (515.68,65) and (518.86,50.55) .. (520,46) .. controls (520.11,45.54) and (519.42,45.21) .. (519,45) .. controls (511.29,41.15) and (511.69,59.76) .. (512,61) .. controls (512.57,63.26) and (516.96,62.13) .. (519,61) .. controls (522.97,58.79) and (525,57.1) .. (525,53) .. controls (525,52.67) and (525.12,53.69) .. (525,54) .. controls (523.55,57.62) and (522.23,61.3) .. (521,65) .. controls (517.31,76.06) and (497.14,64.97) .. (493,60) .. controls (491.49,58.19) and (497.96,59.81) .. (500,61) .. controls (502.88,62.68) and (505.46,64.85) .. (508,67) .. controls (509.8,68.52) and (515.06,70.86) .. (513,72) .. controls (508.92,74.27) and (503.66,72.31) .. (499,72) .. controls (497.63,71.91) and (493.65,70.73) .. (495,71) .. controls (497.64,71.53) and (500.5,71) .. (503,72) .. controls (505.92,73.17) and (513.01,79.42) .. (511,77) .. controls (509.72,75.46) and (500.31,63.96) .. (498,63) .. controls (495.52,61.97) and (490,61.31) .. (490,64) .. controls (490,67.15) and (504.16,71.16) .. (506,73) .. controls (507.84,74.84) and (511,81.6) .. (511,79) .. controls (511,74.15) and (509.17,69.34) .. (507,65) .. controls (506.7,64.4) and (506.89,66.34) .. (507,67) .. controls (507.23,68.36) and (507.52,69.71) .. (508,71) .. controls (508.84,73.23) and (512.13,77.15) .. (510,80) .. controls (509.02,81.31) and (502.15,75.15) .. (503,76) .. controls (504.84,77.84) and (510.71,82.96) .. (509,81) .. controls (504.86,76.27) and (496.7,65.2) .. (489,63) .. controls (487.61,62.6) and (473.16,57.28) .. (477,56) .. controls (481.38,54.54) and (485.06,52.97) .. (489,51) .. controls (489.63,50.69) and (492.87,46.87) .. (494,48) .. controls (495.57,49.57) and (489.49,54.26) .. (489,55) .. controls (488.48,55.78) and (486.06,56.91) .. (487,57) .. controls (495.97,57.82) and (505.94,53.97) .. (514,58) .. controls (515.81,58.91) and (518.54,73.93) .. (518,75) .. controls (516.98,77.03) and (512.68,76.44) .. (511,77) .. controls (510.68,77.11) and (511,78.33) .. (511,78) .. controls (511,71.4) and (515.57,49) .. (525,49) .. controls (529.49,49) and (524.36,57.58) .. (524,59) .. controls (522.55,64.8) and (523.84,71.87) .. (520,77) .. controls (518.25,79.33) and (511.17,88.02) .. (510,81) .. controls (509.14,75.86) and (509.8,73.99) .. (513,70) .. controls (514.47,68.16) and (516.24,66.57) .. (518,65) .. controls (521.14,62.2) and (542.11,42.81) .. (547,42) .. controls (551.97,41.17) and (559.26,40.64) .. (564,41) .. controls (564.22,41.02) and (569,43) .. (569,43) .. controls (569,43) and (567.71,42.24) .. (567,42) .. controls (562.09,40.36) and (556.84,38.21) .. (552,37) .. controls (551.54,36.89) and (550.67,36.33) .. (551,36) .. controls (551.84,35.16) and (564.48,38.39) .. (566,39) .. controls (567.96,39.78) and (574.03,41.55) .. (572,41) .. controls (549.65,34.9) and (528.18,30.46) .. (506,36) .. controls (500.79,37.3) and (492.68,40.88) .. (488,44) .. controls (485.8,45.47) and (484.26,46.87) .. (482,48) .. controls (481.93,48.03) and (478.93,51.03) .. (479,51) .. controls (497.28,41.86) and (519.51,24.76) .. (542,36) .. controls (543.33,36.66) and (537.36,37.91) .. (537,38) .. controls (534.84,38.54) and (533.01,40.5) .. (531,41) .. controls (520.32,43.67) and (508.96,40.97) .. (498,42) .. controls (496.06,42.18) and (493,43.06) .. (493,45) .. controls (493,46.67) and (496.4,45.46) .. (498,45) .. controls (499.77,44.49) and (508.4,40.56) .. (511,39) .. controls (512.21,38.27) and (512.59,36.12) .. (514,36) .. controls (526.44,34.9) and (564.73,35.36) .. (578,42) .. controls (578.3,42.15) and (577.33,42.03) .. (577,42) .. controls (574.1,41.71) and (566.51,42.12) .. (563,43) .. controls (560.63,43.59) and (535.91,51.09) .. (534,53) ;
\draw  [color={rgb, 255:red, 155; green, 155; blue, 155 }  ,draw opacity=1 ][line width=3] [line join = round][line cap = round] (468,59) .. controls (474.71,59) and (492.63,70.63) .. (497,75) .. controls (499.77,77.77) and (504.23,80.23) .. (507,83) .. controls (508,84) and (511,87) .. (510,86) .. controls (507.49,83.49) and (506.83,81.83) .. (505,80) .. controls (500.94,75.94) and (494.07,77.07) .. (491,74) .. controls (489.37,72.37) and (487.48,68.24) .. (485,67) .. controls (483.2,66.1) and (471,60.89) .. (471,60) .. controls (471,58.63) and (473.77,60.39) .. (475,61) .. controls (479.86,63.43) and (485.37,68.37) .. (490,73) .. controls (490.53,73.53) and (489.53,71.53) .. (489,71) .. controls (488.15,70.15) and (486.67,70) .. (486,69) .. controls (484.25,66.38) and (480.18,64) .. (477,64) .. controls (475.2,64) and (480.73,64.73) .. (482,66) .. controls (484.56,68.56) and (485.17,69.11) .. (488,71) .. controls (494.54,75.36) and (501.33,82.33) .. (506,87) .. controls (507.41,88.41) and (510,87.77) .. (510,90) ;
\draw  [color={rgb, 255:red, 155; green, 155; blue, 155 }  ,draw opacity=1 ][line width=3] [line join = round][line cap = round] (591,47) .. controls (592.15,47) and (593.96,48.49) .. (596,49) .. controls (597.37,49.34) and (600,53) .. (599,52) .. controls (597.21,50.21) and (588.14,43) .. (590,43) .. controls (591.21,43) and (596.36,46.91) .. (598,48) .. controls (599,48.67) and (599.15,50.15) .. (600,51) .. controls (600.53,51.53) and (602.57,52.48) .. (602,52) .. controls (599.64,50.03) and (597.46,47.84) .. (595,46) .. controls (594.04,45.28) and (593.17,47.71) .. (592,48) .. controls (590.06,48.49) and (587.97,48.33) .. (586,48) .. controls (580.54,47.09) and (577.17,38.65) .. (561,40) .. controls (559.63,40.11) and (562.61,42.77) .. (562,44) .. controls (561.65,44.71) and (553.23,50.39) .. (552,51) .. controls (538.66,57.67) and (526.83,71.42) .. (521,86) .. controls (519.55,89.62) and (516,91.65) .. (516,95) .. controls (516,96.67) and (519.04,93.37) .. (520,92) .. controls (521.54,89.8) and (522.58,87.28) .. (524,85) .. controls (527.05,80.13) and (525.37,77.5) .. (526,70) .. controls (526.3,66.4) and (529.74,62.9) .. (532,60) .. controls (533.48,58.1) and (534.93,56.15) .. (536,54) .. controls (536.26,53.49) and (536.37,49.21) .. (537,49) .. controls (543.02,46.99) and (550.44,48.09) .. (557,47) .. controls (560.64,46.39) and (566.27,46.25) .. (570,47) .. controls (571.31,47.26) and (573.4,45.81) .. (574,47) .. controls (576.11,51.22) and (572.13,47.98) .. (564,49) .. controls (552.7,50.41) and (543.52,51.48) .. (536,59) ;
\draw  [color={rgb, 255:red, 155; green, 155; blue, 155 }  ,draw opacity=1 ][line width=3] [line join = round][line cap = round] (513,28) .. controls (518.3,24.47) and (525.64,26.35) .. (532,26) .. controls (544.34,25.31) and (553.49,27.5) .. (566,30) .. controls (568.23,30.45) and (574.1,35.55) .. (577,37) .. controls (580.76,38.88) and (593.32,45.63) .. (595,49) .. controls (595.33,49.67) and (593.74,48.11) .. (593,48) .. controls (591.02,47.72) and (588.98,48.28) .. (587,48) .. controls (585.54,47.79) and (579.91,45.91) .. (579,45) ;
\draw  [color={rgb, 255:red, 155; green, 155; blue, 155 }  ,draw opacity=1 ][line width=3] [line join = round][line cap = round] (538,67) .. controls (538,67) and (538,67) .. (538,67) ;
\draw [line width=3]    (523,95) .. controls (524,81) and (549,35) .. (605,58) ;
\draw  [color={rgb, 255:red, 155; green, 155; blue, 155 }  ,draw opacity=1 ][line width=3] [line join = round][line cap = round] (469,59) .. controls (473.79,63.79) and (480.24,63.36) .. (484,69) .. controls (484.26,69.39) and (483.42,68.21) .. (483,68) .. controls (481.03,67.01) and (479.21,67.21) .. (478,66) .. controls (477.76,65.76) and (477,66) .. (477,66) .. controls (477,66) and (484.48,70.58) .. (485,71) .. controls (489.06,74.25) and (503.84,89) .. (507,89) ;
\draw  [color={rgb, 255:red, 155; green, 155; blue, 155 }  ,draw opacity=1 ][line width=3] [line join = round][line cap = round] (476,67) .. controls (473.59,64.59) and (473.31,62.31) .. (471,60) ;
\draw  [color={rgb, 255:red, 155; green, 155; blue, 155 }  ,draw opacity=1 ][line width=3] [line join = round][line cap = round] (466,58) .. controls (464.19,58) and (472.55,63.75) .. (474,65) ;
\draw  [color={rgb, 255:red, 155; green, 155; blue, 155 }  ,draw opacity=1 ][line width=3] [line join = round][line cap = round] (477,66) .. controls (479.54,68.54) and (486.72,76) .. (490,76) ;
\draw  [color={rgb, 255:red, 74; green, 74; blue, 74 }  ,draw opacity=1 ][line width=3] [line join = round][line cap = round] (595,58) .. controls (596.16,58) and (598.34,58.67) .. (599,59) .. controls (599.94,59.47) and (602.94,59.53) .. (602,60) .. controls (599.67,61.17) and (596.29,57.35) .. (595,58) .. controls (593.01,59) and (602.03,61) .. (595,61) ;

\draw (204,131.4) node [anchor=north west][inner sep=0.75pt]  [font=\Large,xscale=0.75,yscale=0.75]  {$\bigodot $};
\draw (419,139.4) node [anchor=north west][inner sep=0.75pt]  [font=\fontsize{4.41em}{5.29em}\selectfont,xscale=0.75,yscale=0.75]  {$=$};
\draw (56,40.4) node [anchor=north west][inner sep=0.75pt]  [xscale=0.75,yscale=0.75]  {$1$};
\draw (187,44.4) node [anchor=north west][inner sep=0.75pt]  [xscale=0.75,yscale=0.75]  {$2$};
\draw (125,88.4) node [anchor=north west][inner sep=0.75pt]  [xscale=0.75,yscale=0.75]  {$3$};
\draw (124,184.4) node [anchor=north west][inner sep=0.75pt]  [xscale=0.75,yscale=0.75]  {$4$};
\draw (237,44.4) node [anchor=north west][inner sep=0.75pt]  [xscale=0.75,yscale=0.75]  {$1$};
\draw (368,48.4) node [anchor=north west][inner sep=0.75pt]  [xscale=0.75,yscale=0.75]  {$2$};
\draw (306,92.4) node [anchor=north west][inner sep=0.75pt]  [xscale=0.75,yscale=0.75]  {$3$};
\draw (307,185.4) node [anchor=north west][inner sep=0.75pt]  [xscale=0.75,yscale=0.75]  {$4$};
\draw (449,45.4) node [anchor=north west][inner sep=0.75pt]  [xscale=0.75,yscale=0.75]  {$1$};
\draw (614,48.4) node [anchor=north west][inner sep=0.75pt]  [xscale=0.75,yscale=0.75]  {$2$};
\draw (518,188.4) node [anchor=north west][inner sep=0.75pt]  [xscale=0.75,yscale=0.75]  {$4$};
\draw (519,94.4) node [anchor=north west][inner sep=0.75pt]  [xscale=0.75,yscale=0.75]  {$3$};

\end{tikzpicture}

%% file: color.tikz
\tikzset{every picture/.style={line width=0.75pt}} 

\begin{tikzpicture}[x=0.75pt,y=0.75pt,yscale=-1,xscale=1]
\hspace{10em}
\draw [line width=2.25]    (137.47,82.73) -- (137.47,180.69) ;
\draw [line width=2.25]    (141.82,78.41) .. controls (178.4,118.51) and (153.14,170.33) .. (144.43,182.42) ;
\draw [color={rgb, 255:red, 208; green, 2; blue, 27 }  ,draw opacity=1 ][line width=2.25]    (131.37,80.14) .. controls (100.89,115.91) and (115.69,167.74) .. (131.37,182.42) ;
\draw [line width=2.25]    (142.69,73.59) -- (277.68,153.05) ;
\draw [line width=2.25]    (145.3,192.42) -- (275.07,160.83) ;

\draw (131.53,62.57) node [anchor=north west][inner sep=0.75pt]   [align=left] {1};
\draw (132.4,183.85) node [anchor=north west][inner sep=0.75pt]   [align=left] {2};
\draw (277.84,149.8) node [anchor=north west][inner sep=0.75pt]   [align=left] {3};

\end{tikzpicture}

%% file: color2.tikz
\tikzset{every picture/.style={line width=0.75pt}} 

\begin{tikzpicture}[x=0.75pt,y=0.75pt,yscale=-1,xscale=1]
\hspace{10em}
\draw [line width=2.25]    (137.47,82.73) -- (137.47,180.69) ;
\draw [line width=2.25]    (141.82,78.41) .. controls (178.4,118.51) and (153.14,170.33) .. (144.43,182.42) ;
\draw [color={rgb, 255:red, 208; green, 2; blue, 27 }  ,draw opacity=1 ][line width=2.25]    (131.37,80.14) .. controls (100.89,115.91) and (115.69,167.74) .. (131.37,182.42) ;
\draw [color={rgb, 255:red, 208; green, 2; blue, 27 }  ,draw opacity=1 ][line width=2.25]    (142.69,73.59) -- (277.68,153.05) ;
\draw [line width=2.25]    (145.3,192.42) -- (275.07,160.83) ;
\draw [line width=2.25]    (146,67) .. controls (207,51.42) and (275,125.42) .. (282,144.42) ;

\draw (131.53,62.57) node [anchor=north west][inner sep=0.75pt]   [align=left] {1};
\draw (132.4,183.85) node [anchor=north west][inner sep=0.75pt]   [align=left] {2};
\draw (277.84,149.8) node [anchor=north west][inner sep=0.75pt]   [align=left] {3};

\end{tikzpicture}

%% file: Firstfiltration.tikz
\tikzset{every picture/.style={line width=0.75pt}} 
\hspace{10em}
\begin{tikzpicture}[x=0.75pt,y=0.75pt,yscale=-1,xscale=1]

\draw [line width=3]    (206,69) -- (274,45) ;
\draw [line width=3]    (285,43) -- (360,64) ;
\draw [line width=3]    (215,125) -- (355,125) ;
\draw [line width=3]    (225.5,214.69) -- (280.47,190.23) ;
\draw [line width=3]    (291.37,187.75) -- (344.65,215.41) ;
\draw [line width=3]    (227.77,220.6) -- (339.79,218.63) ;
\draw  [fill={rgb, 255:red, 155; green, 155; blue, 155 }  ,fill opacity=1 ][dash pattern={on 0.84pt off 2.51pt}] (339.48,215.77) -- (231,216) -- (289.32,191.82) -- cycle ;

\draw (194,60.4) node [anchor=north west][inner sep=0.75pt]  [xscale=0.75,yscale=0.75]  {$1$};
\draw (362,56.4) node [anchor=north west][inner sep=0.75pt]  [xscale=0.75,yscale=0.75]  {$3$};
\draw (275,36.4) node [anchor=north west][inner sep=0.75pt]  [xscale=0.75,yscale=0.75]  {$2$};
\draw (203,113.4) node [anchor=north west][inner sep=0.75pt]  [xscale=0.75,yscale=0.75]  {$1$};
\draw (359,114.4) node [anchor=north west][inner sep=0.75pt]  [xscale=0.75,yscale=0.75]  {$3$};
\draw (215.18,207.54) node [anchor=north west][inner sep=0.75pt]  [rotate=-357.46,xscale=0.75,yscale=0.75]  {$1$};
\draw (345.31,207.78) node [anchor=north west][inner sep=0.75pt]  [rotate=-357.46,xscale=0.75,yscale=0.75]  {$3$};
\draw (281.96,177.64) node [anchor=north west][inner sep=0.75pt]  [rotate=-357.46,xscale=0.75,yscale=0.75]  {$2$};
\draw (408,48.4) node [anchor=north west][inner sep=0.75pt]  [xscale=0.75,yscale=0.75]  {$G$};
\draw (408,107.4) node [anchor=north west][inner sep=0.75pt]  [xscale=0.75,yscale=0.75]  {$H$};
\draw (409,178.4) node [anchor=north west][inner sep=0.75pt]  [xscale=0.75,yscale=0.75]  {$K$};

\end{tikzpicture}

%% file: Secondfiltration.tikz
\tikzset{every picture/.style={line width=0.75pt}} 
\hspace{2em}
\begin{tikzpicture}[x=0.75pt,y=0.75pt,yscale=-1,xscale=1]

\draw [line width=3]    (25,68) -- (80,46) ;
\draw [line width=3]    (91,44) -- (156,67) ;
\draw [line width=3]    (23,74) -- (153,74) ;
\draw [line width=3]    (29.5,168.69) -- (84.47,144.23) ;
\draw [line width=3]    (95.37,141.75) -- (148.65,169.41) ;
\draw [line width=3]    (31.77,174.6) -- (143.79,172.63) ;
\draw  [fill={rgb, 255:red, 155; green, 155; blue, 155 }  ,fill opacity=1 ][dash pattern={on 0.84pt off 2.51pt}] (143.48,169.77) -- (35,170) -- (93.32,145.82) -- cycle ;
\draw  [fill={rgb, 255:red, 155; green, 155; blue, 155 }  ,fill opacity=0.29 ][dash pattern={on 0.84pt off 2.51pt}] (156.24,72.87) -- (28,73) -- (82.81,45.83) -- cycle ;
\draw [line width=3]  [dash pattern={on 3.38pt off 3.27pt}]  (297,16) -- (297,192) ;
\draw [line width=3]    (357,64) -- (412,42) ;
\draw [line width=3]    (423,40) -- (488,63) ;
\draw [line width=3]    (361.5,164.69) -- (416.47,140.23) ;
\draw [line width=3]    (427.37,137.75) -- (480.65,165.41) ;
\draw [line width=3]    (363.77,170.6) -- (475.79,168.63) ;
\draw  [fill={rgb, 255:red, 155; green, 155; blue, 155 }  ,fill opacity=1 ][dash pattern={on 0.84pt off 2.51pt}] (475.48,165.77) -- (367,166) -- (425.32,141.82) -- cycle ;
\draw [line width=3]    (364,176) .. controls (397,204) and (442,204) .. (482,174) ;

\draw (13,62.4) node [anchor=north west][inner sep=0.75pt]  [xscale=0.75,yscale=0.75]  {$1$};
\draw (156,62.4) node [anchor=north west][inner sep=0.75pt]  [xscale=0.75,yscale=0.75]  {$3$};
\draw (81,37.4) node [anchor=north west][inner sep=0.75pt]  [xscale=0.75,yscale=0.75]  {$2$};
\draw (19.18,161.54) node [anchor=north west][inner sep=0.75pt]  [rotate=-357.46,xscale=0.75,yscale=0.75]  {$1$};
\draw (149.31,161.78) node [anchor=north west][inner sep=0.75pt]  [rotate=-357.46,xscale=0.75,yscale=0.75]  {$3$};
\draw (85.96,131.64) node [anchor=north west][inner sep=0.75pt]  [rotate=-357.46,xscale=0.75,yscale=0.75]  {$2$};
\draw (188,48.4) node [anchor=north west][inner sep=0.75pt]  [xscale=0.75,yscale=0.75]  {$G\odot H$};
\draw (203,132.4) node [anchor=north west][inner sep=0.75pt]  [xscale=0.75,yscale=0.75]  {$K$};
\draw (345,58.4) node [anchor=north west][inner sep=0.75pt]  [xscale=0.75,yscale=0.75]  {$1$};
\draw (488,58.4) node [anchor=north west][inner sep=0.75pt]  [xscale=0.75,yscale=0.75]  {$3$};
\draw (413,33.4) node [anchor=north west][inner sep=0.75pt]  [xscale=0.75,yscale=0.75]  {$2$};
\draw (351.18,157.54) node [anchor=north west][inner sep=0.75pt]  [rotate=-357.46,xscale=0.75,yscale=0.75]  {$1$};
\draw (481.31,157.78) node [anchor=north west][inner sep=0.75pt]  [rotate=-357.46,xscale=0.75,yscale=0.75]  {$3$};
\draw (417.96,127.64) node [anchor=north west][inner sep=0.75pt]  [rotate=-357.46,xscale=0.75,yscale=0.75]  {$2$};
\draw (527,46.4) node [anchor=north west][inner sep=0.75pt]  [xscale=0.75,yscale=0.75]  {$G$};
\draw (515,133.4) node [anchor=north west][inner sep=0.75pt]  [xscale=0.75,yscale=0.75]  {$K\odot H$};

\end{tikzpicture}

%% file: Thirdfiltration.tikz
\tikzset{every picture/.style={line width=0.75pt}} 
\hspace{14em}
\begin{tikzpicture}[x=0.75pt,y=0.75pt,yscale=-1,xscale=1]

\draw [line width=3]    (361.5,164.69) -- (416.47,140.23) ;
\draw [line width=3]    (427.37,137.75) -- (480.65,165.41) ;
\draw [line width=3]    (363.77,170.6) -- (475.79,168.63) ;
\draw  [fill={rgb, 255:red, 155; green, 155; blue, 155 }  ,fill opacity=1 ][dash pattern={on 0.84pt off 2.51pt}] (475.48,165.77) -- (367,166) -- (425.32,141.82) -- cycle ;
\draw [line width=3]    (364,176) .. controls (397,204) and (442,204) .. (482,174) ;
\draw  [draw opacity=0][fill={rgb, 255:red, 0; green, 0; blue, 0 }  ,fill opacity=0.18 ][dash pattern={on 0.84pt off 2.51pt}] (475.74,169.28) .. controls (474.57,182.41) and (451.12,192.93) .. (422.35,193) .. controls (392.8,193.08) and (368.82,182.1) .. (368.79,168.48) .. controls (368.79,168.23) and (368.8,167.97) .. (368.81,167.72) -- (422.29,168.35) -- cycle ; \draw  [dash pattern={on 0.84pt off 2.51pt}] (475.74,169.28) .. controls (474.57,182.41) and (451.12,192.93) .. (422.35,193) .. controls (392.8,193.08) and (368.82,182.1) .. (368.79,168.48) .. controls (368.79,168.23) and (368.8,167.97) .. (368.81,167.72) ;  

\draw (351.18,157.54) node [anchor=north west][inner sep=0.75pt]  [rotate=-357.46,xscale=0.75,yscale=0.75]  {$1$};
\draw (481.31,157.78) node [anchor=north west][inner sep=0.75pt]  [rotate=-357.46,xscale=0.75,yscale=0.75]  {$3$};
\draw (417.96,127.64) node [anchor=north west][inner sep=0.75pt]  [rotate=-357.46,xscale=0.75,yscale=0.75]  {$2$};
\draw (515,133.4) node [anchor=north west][inner sep=0.75pt]  [xscale=0.75,yscale=0.75]  {$G\odot H\odot K$};

\end{tikzpicture}

%% file: holes.tikz
\tikzset{every picture/.style={line width=0.75pt}} 

\begin{tikzpicture}[x=0.65pt,y=0.65pt,yscale=-1,xscale=1]

\draw [line width=3]    (23,154) .. controls (23.05,129.03) and (71.62,109.03) .. (95,128) ;
\draw [line width=3]    (28,160) -- (96,136) ;
\draw [line width=3]    (28,169) -- (92,210) ;
\draw [line width=3]    (108,210) -- (171,167) ;
\draw [line width=3]    (115,136) -- (170,155) ;
\draw [line width=3]    (247.83,155) .. controls (247.89,130.03) and (296.45,110.03) .. (319.83,129) ;
\draw [line width=3]    (252.83,161) -- (320.83,137) ;
\draw [line width=3]    (252.83,170) -- (316.83,211) ;
\draw [line width=3]    (335.83,213) -- (395.83,168) ;
\draw [line width=3]    (339.83,137) -- (394.83,156) ;
\draw [line width=3]    (329,214) -- (329,144) ;
\draw [line width=3]    (459.83,164) .. controls (457,138) and (508.62,111.03) .. (532,130) ;
\draw [line width=3]    (464.83,170) -- (533,136) ;
\draw [line width=3]    (464.83,179) -- (528.83,220) ;
\draw [line width=3]    (547.83,222) -- (607.83,177) ;
\draw [line width=3]    (552,135) -- (606.83,165) ;
\draw [line width=3]    (468,175) -- (536,174) ;
\draw [line width=3]    (548,174) -- (605.17,172) ;
\draw  [fill={rgb, 255:red, 155; green, 155; blue, 155 }  ,fill opacity=0.28 ][dash pattern={on 0.84pt off 2.51pt}] (252.83,165.04) -- (325,140) -- (325,210.14) -- cycle ;
\draw  [fill={rgb, 255:red, 128; green, 128; blue, 128 }  ,fill opacity=1 ][dash pattern={on 0.84pt off 2.51pt}] (400,163.13) -- (332,138.25) -- (332,209) -- cycle ;
\draw  [fill={rgb, 255:red, 85; green, 77; blue, 77 }  ,fill opacity=1 ][dash pattern={on 0.84pt off 2.51pt}] (545,138) -- (609,169) -- (545,169) -- cycle ;
\draw  [fill={rgb, 255:red, 85; green, 77; blue, 77 }  ,fill opacity=1 ][dash pattern={on 0.84pt off 2.51pt}] (546,215) -- (607.83,177) -- (546,177) -- cycle ;
\draw  [fill={rgb, 255:red, 85; green, 77; blue, 77 }  ,fill opacity=1 ][dash pattern={on 0.84pt off 2.51pt}] (538,222) -- (472,179) -- (538,179) -- cycle ;
\draw  [fill={rgb, 255:red, 85; green, 77; blue, 77 }  ,fill opacity=1 ][dash pattern={on 0.84pt off 2.51pt}] (535.26,140.32) -- (468.74,171.68) -- (536.66,168.29) -- cycle ;
\draw [line width=3]    (542,226) -- (540,136) ;

\draw (10,105) node [anchor=north west][inner sep=0.75pt]  [font=\large,xscale=0.75,yscale=0.75] [align=left] {\textbf{a)}};
\draw (431,104) node [anchor=north west][inner sep=0.75pt]  [font=\large,xscale=0.75,yscale=0.75] [align=left] {\textbf{c)}};
\draw (221,106) node [anchor=north west][inner sep=0.75pt]  [font=\large,xscale=0.75,yscale=0.75] [align=left] {\textbf{b)}};
\draw (19,159.4) node [anchor=north west][inner sep=0.75pt]  [xscale=0.75,yscale=0.75]  {$1$};
\draw (102,125.4) node [anchor=north west][inner sep=0.75pt]  [xscale=0.75,yscale=0.75]  {$2$};
\draw (175,152.4) node [anchor=north west][inner sep=0.75pt]  [xscale=0.75,yscale=0.75]  {$3$};
\draw (94,207.4) node [anchor=north west][inner sep=0.75pt]  [xscale=0.75,yscale=0.75]  {$4$};
\draw (324,122.4) node [anchor=north west][inner sep=0.75pt]  [xscale=0.75,yscale=0.75]  {$2$};
\draw (324,218.54) node [anchor=north west][inner sep=0.75pt]  [xscale=0.75,yscale=0.75]  {$4$};
\draw (242,156.4) node [anchor=north west][inner sep=0.75pt]  [xscale=0.75,yscale=0.75]  {$1$};
\draw (534,121.4) node [anchor=north west][inner sep=0.75pt]  [xscale=0.75,yscale=0.75]  {$2$};
\draw (536,228.54) node [anchor=north west][inner sep=0.75pt]  [xscale=0.75,yscale=0.75]  {$4$};
\draw (454,164.4) node [anchor=north west][inner sep=0.75pt]  [xscale=0.75,yscale=0.75]  {$1$};
\draw (611.83,163.4) node [anchor=north west][inner sep=0.75pt]  [xscale=0.75,yscale=0.75]  {$3$};
\draw (399,153.4) node [anchor=north west][inner sep=0.75pt]  [xscale=0.75,yscale=0.75]  {$3$};

\end{tikzpicture}

%% file: Fourthfiltration.tikz
\tikzset{every picture/.style={line width=0.75pt}} 
\hspace{2em}
\begin{tikzpicture}[x=0.75pt,y=0.75pt,yscale=-1,xscale=1]

\draw [line width=3]    (165.17,122.5) .. controls (181.33,140) and (206.17,147.5) .. (245.17,140.5) ;
\draw [line width=3]    (258,139) .. controls (278.17,135.5) and (293.17,126.5) .. (305.17,117.5) ;
\draw [line width=3]    (162,56) .. controls (187.17,41.5) and (213.17,38.5) .. (227.17,41.5) ;
\draw [line width=3]    (242.17,41.5) .. controls (258.17,40.5) and (273.17,42.5) .. (289.17,48.5) ;
\draw [line width=3]    (233.17,48.5) .. controls (232.33,58) and (239.17,91.5) .. (253.17,94.5) ;
\draw [line width=3]    (219,197) .. controls (235.17,201.5) and (270.17,195.5) .. (283.17,187.5) ;
\draw [line width=3]    (475,148) .. controls (491.17,152.5) and (526.17,146.5) .. (539.17,138.5) ;
\draw [line width=1.5]  [dash pattern={on 1.69pt off 2.76pt}]  (341,13) -- (341.17,214.5) ;
\draw [line width=3]    (162,303) .. controls (187.17,288.5) and (213.17,285.5) .. (227.17,288.5) ;
\draw [line width=3]    (242.17,288.5) .. controls (258.17,287.5) and (273.17,289.5) .. (289.17,295.5) ;
\draw [line width=3]    (233.17,295.5) .. controls (232.33,305) and (239.17,338.5) .. (253.17,341.5) ;
\draw [line width=1.5]  [dash pattern={on 1.69pt off 2.76pt}]  (341,260) -- (341.17,461.5) ;
\draw  [color={rgb, 255:red, 255; green, 255; blue, 255 }  ,draw opacity=1 ][line width=5.25] [line join = round][line cap = round] (321.17,436.5) .. controls (320.2,434.58) and (316.44,430.89) .. (316.17,429.5) .. controls (315.52,426.24) and (317.56,423.76) .. (314.17,421.5) .. controls (312.91,420.66) and (313.18,423.49) .. (313.17,423.5) .. controls (308.51,425.83) and (305.1,430.67) .. (303.17,435.5) .. controls (302.67,436.74) and (303.17,438.17) .. (303.17,439.5) .. controls (303.17,439.83) and (302.89,438.69) .. (303.17,438.5) .. controls (306.35,436.27) and (309.93,434.66) .. (313.17,432.5) .. controls (316.24,430.45) and (316.52,425.96) .. (318.17,423.5) .. controls (318.69,422.72) and (319.5,422.17) .. (320.17,421.5) .. controls (320.5,421.17) and (319.23,422.03) .. (319.17,422.5) .. controls (318.88,424.82) and (318.88,427.18) .. (319.17,429.5) .. controls (319.6,433.01) and (327.02,426.71) .. (327.17,426.5) .. controls (329.17,423.5) and (320.17,429.5) .. (318.17,432.5) .. controls (315.39,436.67) and (308.17,433.83) .. (303.17,433.5) .. controls (302.23,433.44) and (304.25,431.73) .. (305.17,431.5) .. controls (307.49,430.92) and (308.24,431.14) .. (310.17,430.5) .. controls (310.87,430.26) and (312.91,429.56) .. (312.17,429.5) .. controls (309.84,429.31) and (295.17,432.64) .. (295.17,424.5) ;
\draw  [color={rgb, 255:red, 255; green, 255; blue, 255 }  ,draw opacity=1 ][line width=5.25] [line join = round][line cap = round] (306.17,409.5) .. controls (302.22,409.5) and (300.21,425.16) .. (301.17,428.5) .. controls (301.35,429.14) and (302.55,428.75) .. (303.17,428.5) .. controls (304.92,427.8) and (311.09,417.48) .. (307.17,416.5) .. controls (300.43,414.82) and (296.43,424.63) .. (302.17,427.5) .. controls (305.8,429.32) and (310.13,417.81) .. (308.17,416.5) .. controls (306.46,415.36) and (299.99,422.23) .. (299.17,422.5) .. controls (298.51,422.72) and (296.83,423.28) .. (296.17,423.5) .. controls (295.53,423.71) and (293.5,423.5) .. (294.17,423.5) .. controls (295.48,423.5) and (306.17,418.04) .. (306.17,420.5) .. controls (306.17,421.25) and (298.19,424.43) .. (298.17,424.5) .. controls (295.98,429.96) and (299.05,434.5) .. (304.17,434.5) ;
\draw [line width=1.5]  [dash pattern={on 1.69pt off 2.76pt}]  (517.17,240.5) -- (140.17,241.5) ;
\draw [line width=3]    (434.87,63.82) -- (494.05,43.69) ;
\draw [line width=3]    (434.87,71.37) -- (494,109) ;
\draw [line width=3]    (507.1,107.43) -- (559.32,69.69) ;
\draw [line width=3]    (510.58,43.69) -- (558.45,59.63) ;
\draw [line width=3]    (501,111) -- (501,50) ;
\draw  [fill={rgb, 255:red, 155; green, 155; blue, 155 }  ,fill opacity=0.28 ][dash pattern={on 0.84pt off 2.51pt}] (434.87,67.21) -- (497.68,46.21) -- (497.68,105.04) -- cycle ;
\draw  [fill={rgb, 255:red, 128; green, 128; blue, 128 }  ,fill opacity=1 ][dash pattern={on 0.84pt off 2.51pt}] (562.95,65.6) -- (503.77,44.75) -- (503.77,104.08) -- cycle ;
\draw [line width=3]    (426.49,320.44) -- (484.41,294.27) ;
\draw [line width=3]    (426.49,327.37) -- (485.55,366.45) ;
\draw [line width=3]    (497.02,360.47) -- (548.01,325.83) ;
\draw [line width=3]    (496,291) -- (547.16,316.59) ;
\draw [line width=3]    (429.18,324.29) -- (486.96,323.52) ;
\draw [line width=3]    (497.16,323.52) -- (545.74,321.98) ;
\draw  [fill={rgb, 255:red, 85; green, 77; blue, 77 }  ,fill opacity=1 ][dash pattern={on 0.84pt off 2.51pt}] (494.61,295.81) -- (549,319.67) -- (494.61,319.67) -- cycle ;
\draw  [fill={rgb, 255:red, 85; green, 77; blue, 77 }  ,fill opacity=1 ][dash pattern={on 0.84pt off 2.51pt}] (495.46,355.08) -- (548.01,325.83) -- (495.46,325.83) -- cycle ;
\draw  [fill={rgb, 255:red, 85; green, 77; blue, 77 }  ,fill opacity=1 ][dash pattern={on 0.84pt off 2.51pt}] (488.66,360.47) -- (432.58,327.37) -- (488.66,327.37) -- cycle ;
\draw  [fill={rgb, 255:red, 85; green, 77; blue, 77 }  ,fill opacity=1 ][dash pattern={on 0.84pt off 2.51pt}] (486.33,297.58) -- (429.81,321.74) -- (487.52,319.12) -- cycle ;
\draw [line width=3]    (492.37,365.37) -- (490.09,292.41) ;
\draw [line width=3]    (165.87,403.82) -- (230,388) ;
\draw [line width=3]    (167.87,409.37) -- (289,427) ;
\draw [line width=3]    (244.58,388.69) -- (292,421) ;
\draw  [fill={rgb, 255:red, 155; green, 155; blue, 155 }  ,fill opacity=0.28 ][dash pattern={on 0.84pt off 2.51pt}] (232.46,386.79) -- (289.86,425.23) -- (173.57,406.96) -- cycle ;

\draw (228.68,32.72) node [anchor=north west][inner sep=0.75pt]  [font=\normalsize,xscale=0.75,yscale=0.75] [align=left] {2};
\draw (291.95,44.96) node [anchor=north west][inner sep=0.75pt]  [font=\normalsize,xscale=0.75,yscale=0.75] [align=left] {3};
\draw (151.44,49.82) node [anchor=north west][inner sep=0.75pt]  [font=\normalsize,xscale=0.75,yscale=0.75] [align=left] {1};
\draw (251.33,92.24) node [anchor=north west][inner sep=0.75pt]  [font=\normalsize,xscale=0.75,yscale=0.75] [align=left] {4};
\draw (246.69,132.53) node [anchor=north west][inner sep=0.75pt]  [font=\normalsize,xscale=0.75,yscale=0.75] [align=left] {4};
\draw (303.56,107.14) node [anchor=north west][inner sep=0.75pt]  [font=\normalsize,xscale=0.75,yscale=0.75] [align=left] {3};
\draw (153.98,116.35) node [anchor=north west][inner sep=0.75pt]  [font=\normalsize,xscale=0.75,yscale=0.75] [align=left] {1};
\draw (282.26,182.28) node [anchor=north west][inner sep=0.75pt]  [font=\normalsize,xscale=0.75,yscale=0.75] [align=left] {3};
\draw (209.31,194.38) node [anchor=north west][inner sep=0.75pt]  [font=\normalsize,xscale=0.75,yscale=0.75] [align=left] {1};
\draw (103.52,22.32) node [anchor=north west][inner sep=0.75pt]  [font=\normalsize,xscale=0.75,yscale=0.75]  {$\delta =0$};
\draw (359.07,25.94) node [anchor=north west][inner sep=0.75pt]  [font=\normalsize,xscale=0.75,yscale=0.75]  {$\delta =1$};
\draw (539.26,133.28) node [anchor=north west][inner sep=0.75pt]  [font=\normalsize,xscale=0.75,yscale=0.75] [align=left] {3};
\draw (464.31,144.38) node [anchor=north west][inner sep=0.75pt]  [font=\normalsize,xscale=0.75,yscale=0.75] [align=left] {1};
\draw (228.68,281.72) node [anchor=north west][inner sep=0.75pt]  [font=\normalsize,xscale=0.75,yscale=0.75] [align=left] {2};
\draw (291.95,291.96) node [anchor=north west][inner sep=0.75pt]  [font=\normalsize,xscale=0.75,yscale=0.75] [align=left] {3};
\draw (151.44,296.82) node [anchor=north west][inner sep=0.75pt]  [font=\normalsize,xscale=0.75,yscale=0.75] [align=left] {1};
\draw (251.33,339.24) node [anchor=north west][inner sep=0.75pt]  [font=\normalsize,xscale=0.75,yscale=0.75] [align=left] {4};
\draw (233.69,374.53) node [anchor=north west][inner sep=0.75pt]  [font=\normalsize,xscale=0.75,yscale=0.75] [align=left] {4};
\draw (292.56,418.14) node [anchor=north west][inner sep=0.75pt]  [font=\normalsize,xscale=0.75,yscale=0.75] [align=left] {3};
\draw (152.98,397.35) node [anchor=north west][inner sep=0.75pt]  [font=\normalsize,xscale=0.75,yscale=0.75] [align=left] {1};
\draw (103.52,269.32) node [anchor=north west][inner sep=0.75pt]  [font=\normalsize,xscale=0.75,yscale=0.75]  {$\delta =2$};
\draw (359.07,272.94) node [anchor=north west][inner sep=0.75pt]  [font=\normalsize,xscale=0.75,yscale=0.75]  {$\delta =3$};
\draw (417.31,57.38) node [anchor=north west][inner sep=0.75pt]  [font=\normalsize,xscale=0.75,yscale=0.75] [align=left] {1};
\draw (495.68,29.72) node [anchor=north west][inner sep=0.75pt]  [font=\normalsize,xscale=0.75,yscale=0.75] [align=left] {2};
\draw (562.95,56.96) node [anchor=north west][inner sep=0.75pt]  [font=\normalsize,xscale=0.75,yscale=0.75] [align=left] {3};
\draw (492.57,108.75) node [anchor=north west][inner sep=0.75pt]  [font=\normalsize,xscale=0.75,yscale=0.75] [align=left] {4};
\draw (483.67,274.49) node [anchor=north west][inner sep=0.75pt]  [font=\normalsize,xscale=0.75,yscale=0.75] [align=left] {2};
\draw (551.06,314.45) node [anchor=north west][inner sep=0.75pt]  [font=\normalsize,xscale=0.75,yscale=0.75] [align=left] {3};
\draw (412.89,315.37) node [anchor=north west][inner sep=0.75pt]  [font=\normalsize,xscale=0.75,yscale=0.75] [align=left] {1};
\draw (486.01,363.23) node [anchor=north west][inner sep=0.75pt]  [font=\normalsize,xscale=0.75,yscale=0.75] [align=left] {4};

\end{tikzpicture}

%% file: betti1.tikz
\tikzset{every picture/.style={line width=0.75pt}} 
\hspace{4em}
\begin{tikzpicture}[x=0.75pt,y=0.75pt,yscale=-1,xscale=1]

\draw [line width=3]    (131.63,143.95) .. controls (140.5,173) and (194.5,180) .. (213.5,144) ;
\draw [line width=3]    (135,137) -- (208.5,137) ;
\draw [line width=3]    (297.5,136) -- (369.5,83) ;
\draw [line width=3]    (432.5,129) -- (504.5,76) ;
\draw [line width=3]    (431.63,138.95) .. controls (440.5,168) and (494.5,175) .. (513.5,139) ;
\draw [line width=3]    (435,132) -- (508.5,132) ;

\draw (123,131.4) node [anchor=north west][inner sep=0.75pt]  [xscale=0.75,yscale=0.75]  {$1$};
\draw (210,127.4) node [anchor=north west][inner sep=0.75pt]  [xscale=0.75,yscale=0.75]  {$2$};
\draw (290,129.4) node [anchor=north west][inner sep=0.75pt]  [xscale=0.75,yscale=0.75]  {$1$};
\draw (368,73.4) node [anchor=north west][inner sep=0.75pt]  [xscale=0.75,yscale=0.75]  {$3$};
\draw (245,123.4) node [anchor=north west][inner sep=0.75pt]  [font=\Large,xscale=0.75,yscale=0.75]  {$\odot $};
\draw (117,90) node [anchor=north west][inner sep=0.75pt]  [font=\normalsize,xscale=0.75,yscale=0.75] [align=left] {G};
\draw (279,91) node [anchor=north west][inner sep=0.75pt]  [font=\normalsize,xscale=0.75,yscale=0.75] [align=left] {H};
\draw (425,124.4) node [anchor=north west][inner sep=0.75pt]  [xscale=0.75,yscale=0.75]  {$1$};
\draw (502,68.4) node [anchor=north west][inner sep=0.75pt]  [xscale=0.75,yscale=0.75]  {$3$};
\draw (366,127.4) node [anchor=north west][inner sep=0.75pt]  [font=\Large,xscale=0.75,yscale=0.75]  {$=$};
\draw (422,90.4) node [anchor=north west][inner sep=0.75pt]  [font=\normalsize,xscale=0.75,yscale=0.75]  {$G\odot H$};
\draw (509,125.4) node [anchor=north west][inner sep=0.75pt]  [xscale=0.75,yscale=0.75]  {$2$};

\end{tikzpicture}

%% file: betti2.tikz
\tikzset{every picture/.style={line width=0.75pt}} 
\hspace{4em}
   \begin{tikzpicture}[x=0.75pt,y=0.75pt,yscale=-1,xscale=1]

\draw [line width=3]    (131.63,143.95) .. controls (140.5,173) and (194.5,180) .. (213.5,144) ;
\draw [line width=3]    (135,137) -- (208.5,137) ;
\draw [line width=3]    (297.5,136) -- (340.5,102) ;
\draw [line width=3]    (451.63,140.95) .. controls (460.5,170) and (514.5,177) .. (533.5,141) ;
\draw [line width=3]    (455,134) -- (528.5,134) ;
\draw [line width=3]    (386.5,135) -- (348.5,103) ;
\draw [line width=3]    (454.5,125) -- (490.5,92) ;
\draw [line width=3]    (531.5,128) -- (498.5,93) ;
\draw  [fill={rgb, 255:red, 155; green, 155; blue, 155 }  ,fill opacity=1 ][dash pattern={on 0.84pt off 2.51pt}] (493.25,94) -- (532,130) -- (454.5,130) -- cycle ;

\draw (123,131.4) node [anchor=north west][inner sep=0.75pt]  [xscale=0.75,yscale=0.75]  {$1$};
\draw (210,127.4) node [anchor=north west][inner sep=0.75pt]  [xscale=0.75,yscale=0.75]  {$2$};
\draw (290,129.4) node [anchor=north west][inner sep=0.75pt]  [xscale=0.75,yscale=0.75]  {$1$};
\draw (338,89.4) node [anchor=north west][inner sep=0.75pt]  [xscale=0.75,yscale=0.75]  {$3$};
\draw (245,123.4) node [anchor=north west][inner sep=0.75pt]  [font=\Large,xscale=0.75,yscale=0.75]  {$\odot $};
\draw (445,126.4) node [anchor=north west][inner sep=0.75pt]  [xscale=0.75,yscale=0.75]  {$1$};
\draw (403,126.4) node [anchor=north west][inner sep=0.75pt]  [font=\Large,xscale=0.75,yscale=0.75]  {$=$};
\draw (420,91.4) node [anchor=north west][inner sep=0.75pt]  [font=\normalsize,xscale=0.75,yscale=0.75]  {$G\odot H$};
\draw (529,127.4) node [anchor=north west][inner sep=0.75pt]  [xscale=0.75,yscale=0.75]  {$2$};
\draw (384,130.4) node [anchor=north west][inner sep=0.75pt]  [xscale=0.75,yscale=0.75]  {$2$};
\draw (488,79.4) node [anchor=north west][inner sep=0.75pt]  [xscale=0.75,yscale=0.75]  {$3$};
\draw (142,90.4) node [anchor=north west][inner sep=0.75pt]  [xscale=0.75,yscale=0.75]  {$G$};
\draw (289,88.4) node [anchor=north west][inner sep=0.75pt]  [xscale=0.75,yscale=0.75]  {$H$};

\end{tikzpicture}

%% file: betti3.tikz
\tikzset{every picture/.style={line width=0.75pt}} 
\hspace{4em}
\begin{tikzpicture}[x=0.75pt,y=0.75pt,yscale=-1,xscale=1]

\draw [line width=3]    (131.63,143.95) .. controls (140.5,173) and (194.5,180) .. (213.5,144) ;
\draw [line width=3]    (135,137) -- (208.5,137) ;
\draw [line width=3]    (296.5,145) -- (367,89) ;
\draw [line width=3]    (372,165) -- (373,91) ;
\draw [line width=3]    (327,179) -- (367,89) ;
\draw [line width=3]    (299,150) -- (330,158) ;
\draw [line width=3]    (340,162) -- (366,169) ;
\draw [line width=3]    (330,185) -- (367,174) ;
\draw [line width=3]    (445.5,147) -- (516,91) ;
\draw [line width=3]    (521,167) -- (522,93) ;
\draw [line width=3]    (476,181) -- (516,91) ;
\draw [line width=3]    (448,152) -- (479,160) ;
\draw [line width=3]    (489,164) -- (515,171) ;
\draw [line width=3]    (479,187) -- (516,176) ;
\draw [line width=3]    (443,155) -- (469,184) ;
\draw [line width=3]    (436.63,157.95) .. controls (417,174) and (450,208) .. (466,193) ;

\draw (123,131.4) node [anchor=north west][inner sep=0.75pt]  [xscale=0.75,yscale=0.75]  {$1$};
\draw (210,127.4) node [anchor=north west][inner sep=0.75pt]  [xscale=0.75,yscale=0.75]  {$2$};
\draw (286,140.4) node [anchor=north west][inner sep=0.75pt]  [xscale=0.75,yscale=0.75]  {$1$};
\draw (368,73.4) node [anchor=north west][inner sep=0.75pt]  [xscale=0.75,yscale=0.75]  {$3$};
\draw (245,123.4) node [anchor=north west][inner sep=0.75pt]  [font=\Large,xscale=0.75,yscale=0.75]  {$\odot $};
\draw (117,90) node [anchor=north west][inner sep=0.75pt]  [font=\normalsize,xscale=0.75,yscale=0.75] [align=left] {G};
\draw (279,91) node [anchor=north west][inner sep=0.75pt]  [font=\normalsize,xscale=0.75,yscale=0.75] [align=left] {H};
\draw (403,126.4) node [anchor=north west][inner sep=0.75pt]  [font=\Large,xscale=0.75,yscale=0.75]  {$=$};
\draw (428,93) node [anchor=north west][inner sep=0.75pt]  [font=\normalsize,xscale=0.75,yscale=0.75] [align=left] {$G\odot H$};
\draw (366,161.4) node [anchor=north west][inner sep=0.75pt]  [xscale=0.75,yscale=0.75]  {$4$};
\draw (318,176.4) node [anchor=north west][inner sep=0.75pt]  [xscale=0.75,yscale=0.75]  {$2$};
\draw (435,142.4) node [anchor=north west][inner sep=0.75pt]  [xscale=0.75,yscale=0.75]  {$1$};
\draw (517,75.4) node [anchor=north west][inner sep=0.75pt]  [xscale=0.75,yscale=0.75]  {$3$};
\draw (515,163.4) node [anchor=north west][inner sep=0.75pt]  [xscale=0.75,yscale=0.75]  {$4$};
\draw (467,178.4) node [anchor=north west][inner sep=0.75pt]  [xscale=0.75,yscale=0.75]  {$2$};

\end{tikzpicture}